\catcode`\@=11

\font\tenmsa=msam10
\font\sevenmsa=msam7
\font\fivemsa=msam5
\font\tenmsb=msbm10
\font\sevenmsb=msbm7
\font\fivemsb=msbm5
\newfam\msafam
\newfam\msbfam
\textfont\msafam=\tenmsa  \scriptfont\msafam=\sevenmsa
  \scriptscriptfont\msafam=\fivemsa
\textfont\msbfam=\tenmsb  \scriptfont\msbfam=\sevenmsb
  \scriptscriptfont\msbfam=\fivemsb

\def\hexnumber@#1{\ifnum#1<10 \number#1\else
 \ifnum#1=10 A\else\ifnum#1=11 B\else\ifnum#1=12 C\else
 \ifnum#1=13 D\else\ifnum#1=14 E\else\ifnum#1=15 F\fi\fi\fi\fi\fi\fi\fi}

\def\msa@{\hexnumber@\msafam}
\def\msb@{\hexnumber@\msbfam}
\mathchardef\boxdot="2\msa@00
\mathchardef\boxplus="2\msa@01
\mathchardef\boxtimes="2\msa@02
\mathchardef\square="0\msa@03
\mathchardef\blacksquare="0\msa@04
\mathchardef\centerdot="2\msa@05
\mathchardef\lozenge="0\msa@06
\mathchardef\blacklozenge="0\msa@07
\mathchardef\circlearrowright="3\msa@08
\mathchardef\circlearrowleft="3\msa@09
\mathchardef\rightleftharpoons="3\msa@0A
\mathchardef\leftrightharpoons="3\msa@0B
\mathchardef\boxminus="2\msa@0C
\mathchardef\Vdash="3\msa@0D
\mathchardef\Vvdash="3\msa@0E
\mathchardef\vDash="3\msa@0F
\mathchardef\twoheadrightarrow="3\msa@10
\mathchardef\twoheadleftarrow="3\msa@11
\mathchardef\leftleftarrows="3\msa@12
\mathchardef\rightrightarrows="3\msa@13
\mathchardef\upuparrows="3\msa@14
\mathchardef\downdownarrows="3\msa@15
\mathchardef\upharpoonright="3\msa@16

\mathchardef\downharpoonright="3\msa@17
\mathchardef\upharpoonleft="3\msa@18
\mathchardef\downharpoonleft="3\msa@19
\mathchardef\rightarrowtail="3\msa@1A
\mathchardef\leftarrowtail="3\msa@1B
\mathchardef\leftrightarrows="3\msa@1C
\mathchardef\rightleftarrows="3\msa@1D
\mathchardef\Lsh="3\msa@1E
\mathchardef\Rsh="3\msa@1F
\mathchardef\rightsquigarrow="3\msa@20
\mathchardef\leftrightsquigarrow="3\msa@21
\mathchardef\looparrowleft="3\msa@22
\mathchardef\looparrowright="3\msa@23
\mathchardef\circeq="3\msa@24
\mathchardef\succsim="3\msa@25
\mathchardef\gtrsim="3\msa@26
\mathchardef\gtrapprox="3\msa@27
\mathchardef\multimap="3\msa@28
\mathchardef\therefore="3\msa@29
\mathchardef\because="3\msa@2A
\mathchardef\doteqdot="3\msa@2B

\mathchardef\triangleq="3\msa@2C
\mathchardef\precsim="3\msa@2D
\mathchardef\lesssim="3\msa@2E
\mathchardef\lessapprox="3\msa@2F
\mathchardef\eqslantless="3\msa@30
\mathchardef\eqslantgtr="3\msa@31
\mathchardef\curlyeqprec="3\msa@32
\mathchardef\curlyeqsucc="3\msa@33
\mathchardef\preccurlyeq="3\msa@34
\mathchardef\leqq="3\msa@35
\mathchardef\leqslant="3\msa@36
\mathchardef\lessgtr="3\msa@37
\mathchardef\backprime="0\msa@38
\mathchardef\risingdotseq="3\msa@3A
\mathchardef\fallingdotseq="3\msa@3B
\mathchardef\succcurlyeq="3\msa@3C
\mathchardef\geqq="3\msa@3D
\mathchardef\geqslant="3\msa@3E
\mathchardef\gtrless="3\msa@3F
\mathchardef\sqsubset="3\msa@40
\mathchardef\sqsupset="3\msa@41
\mathchardef\trianglerighteq="3\msa@44
\mathchardef\trianglelefteq="3\msa@45
\mathchardef\bigstar="0\msa@46
\mathchardef\between="3\msa@47
\mathchardef\blacktriangledown="0\msa@48
\mathchardef\blacktriangleright="3\msa@49
\mathchardef\blacktriangleleft="3\msa@4A
\mathchardef\blacktriangle="0\msa@4E
\mathchardef\triangledown="0\msa@4F
\mathchardef\eqcirc="3\msa@50
\mathchardef\lesseqgtr="3\msa@51
\mathchardef\gtreqless="3\msa@52
\mathchardef\lesseqqgtr="3\msa@53
\mathchardef\gtreqqless="3\msa@54
\mathchardef\Rrightarrow="3\msa@56
\mathchardef\Lleftarrow="3\msa@57
\mathchardef\veebar="2\msa@59
\mathchardef\barwedge="2\msa@5A
\mathchardef\doublebarwedge="2\msa@5B
\mathchardef\angle="0\msa@5C
\mathchardef\measuredangle="0\msa@5D
\mathchardef\sphericalangle="0\msa@5E
\mathchardef\varpropto="3\msa@5F
\mathchardef\smallsmile="3\msa@60
\mathchardef\smallfrown="3\msa@61
\mathchardef\Subset="3\msa@62
\mathchardef\Supset="3\msa@63
\mathchardef\Cup="2\msa@64

\mathchardef\Cap="2\msa@65

\mathchardef\curlywedge="2\msa@66
\mathchardef\curlyvee="2\msa@67
\mathchardef\leftthreetimes="2\msa@68
\mathchardef\rightthreetimes="2\msa@69
\mathchardef\subseteqq="3\msa@6A
\mathchardef\supseteqq="3\msa@6B
\mathchardef\bumpeq="3\msa@6C
\mathchardef\Bumpeq="3\msa@6D
\mathchardef\lll="3\msa@6E

\mathchardef\ggg="3\msa@6F

\mathchardef\circledS="0\msa@73
\mathchardef\pitchfork="3\msa@74
\mathchardef\dotplus="2\msa@75
\mathchardef\backsim="3\msa@76
\mathchardef\backsimeq="3\msa@77
\mathchardef\complement="0\msa@7B
\mathchardef\intercal="2\msa@7C
\mathchardef\circledcirc="2\msa@7D
\mathchardef\circledast="2\msa@7E
\mathchardef\circleddash="2\msa@7F
\def\ulcorner{\delimiter"4\msa@70\msa@70 }
\def\urcorner{\delimiter"5\msa@71\msa@71 }
\def\llcorner{\delimiter"4\msa@78\msa@78 }
\def\lrcorner{\delimiter"5\msa@79\msa@79 }
\def\yen{\mathhexbox\msa@55 }
\def\checkmark{\mathhexbox\msa@58 }
\def\circledR{\mathhexbox\msa@72 }
\def\maltese{\mathhexbox\msa@7A }
\mathchardef\lvertneqq="3\msb@00
\mathchardef\gvertneqq="3\msb@01
\mathchardef\nleq="3\msb@02
\mathchardef\ngeq="3\msb@03
\mathchardef\nless="3\msb@04
\mathchardef\ngtr="3\msb@05
\mathchardef\nprec="3\msb@06
\mathchardef\nsucc="3\msb@07
\mathchardef\lneqq="3\msb@08
\mathchardef\gneqq="3\msb@09
\mathchardef\nleqslant="3\msb@0A
\mathchardef\ngeqslant="3\msb@0B
\mathchardef\lneq="3\msb@0C
\mathchardef\gneq="3\msb@0D
\mathchardef\npreceq="3\msb@0E
\mathchardef\nsucceq="3\msb@0F
\mathchardef\precnsim="3\msb@10
\mathchardef\succnsim="3\msb@11
\mathchardef\lnsim="3\msb@12
\mathchardef\gnsim="3\msb@13
\mathchardef\nleqq="3\msb@14
\mathchardef\ngeqq="3\msb@15
\mathchardef\precneqq="3\msb@16
\mathchardef\succneqq="3\msb@17
\mathchardef\precnapprox="3\msb@18
\mathchardef\succnapprox="3\msb@19
\mathchardef\lnapprox="3\msb@1A
\mathchardef\gnapprox="3\msb@1B
\mathchardef\nsim="3\msb@1C
\mathchardef\napprox="3\msb@1D
\mathchardef\nsubseteqq="3\msb@22
\mathchardef\nsupseteqq="3\msb@23
\mathchardef\subsetneqq="3\msb@24
\mathchardef\supsetneqq="3\msb@25
\mathchardef\subsetneq="3\msb@28
\mathchardef\supsetneq="3\msb@29
\mathchardef\nsubseteq="3\msb@2A
\mathchardef\nsupseteq="3\msb@2B
\mathchardef\nparallel="3\msb@2C
\mathchardef\nmid="3\msb@2D
\mathchardef\nshortmid="3\msb@2E
\mathchardef\nshortparallel="3\msb@2F
\mathchardef\nvdash="3\msb@30
\mathchardef\nVdash="3\msb@31
\mathchardef\nvDash="3\msb@32
\mathchardef\nVDash="3\msb@33
\mathchardef\ntrianglerighteq="3\msb@34
\mathchardef\ntrianglelefteq="3\msb@35
\mathchardef\ntriangleleft="3\msb@36
\mathchardef\ntriangleright="3\msb@37
\mathchardef\nleftarrow="3\msb@38
\mathchardef\nrightarrow="3\msb@39
\mathchardef\nLeftarrow="3\msb@3A
\mathchardef\nRightarrow="3\msb@3B
\mathchardef\nLeftrightarrow="3\msb@3C
\mathchardef\nleftrightarrow="3\msb@3D
\mathchardef\divideontimes="2\msb@3E
\mathchardef\varnothing="0\msb@3F
\mathchardef\nexists="0\msb@40
\mathchardef\mho="0\msb@66
\mathchardef\thorn="0\msb@67
\mathchardef\beth="0\msb@69
\mathchardef\gimel="0\msb@6A
\mathchardef\daleth="0\msb@6B
\mathchardef\lessdot="3\msb@6C
\mathchardef\gtrdot="3\msb@6D
\mathchardef\ltimes="2\msb@6E
\mathchardef\rtimes="2\msb@6F
\mathchardef\shortmid="3\msb@70
\mathchardef\shortparallel="3\msb@71
\mathchardef\smallsetminus="2\msb@72
\mathchardef\thicksim="3\msb@73
\mathchardef\thickapprox="3\msb@74
\mathchardef\approxeq="3\msb@75
\mathchardef\succapprox="3\msb@76
\mathchardef\precapprox="3\msb@77
\mathchardef\curvearrowleft="3\msb@78
\mathchardef\curvearrowright="3\msb@79
\mathchardef\digamma="0\msb@7A
\mathchardef\varkappa="0\msb@7B
\mathchardef\hslash="0\msb@7D
\mathchardef\hbar="0\msb@7E
\mathchardef\backepsilon="3\msb@7F
\def\Bbb{\ifmmode\let\next\Bbb@\else
 \def\next{\errmessage{Use \string\Bbb\space only in math mode}}\fi\next}
\def\Bbb@#1{{\Bbb@@{#1}}}
\def\Bbb@@#1{\fam\msbfam#1}

\catcode`\@=\active

\def\del{\partial}

 \def\CC{\hbox{{$\cal C$}}}
 
\def\CL{\hbox{{$\cal L$}}}

\def\CR{\hbox{{$\cal R$}}} 
 \def\CV{\hbox{{$\cal V$}}}
\def\CM{\hbox{{$\cal M$}}} 
  
\def\CO{\hbox{{$\cal O$}}} 
\def\CX{\hbox{{$\cal X$}}} \def\CZ{\hbox{{$\cal Z$}}}

\def\lform{\hbox{$\sqcup$}\llap{\hbox{$\sqcap$}}}

\def\C{{\Bbb C}}
\def\Z{{\Bbb Z}}

\def\eps{{\epsilon}}

\def\trace{{\rm Tr\, }}
\def\aut{{\rm Aut\, }}

\def\ant{{{\scriptstyle S}}}

\def\cross{{\triangleright\!\!\!<}}
\def\cocross{{>\!\!\!\triangleleft}}

\def\rbiprod{{\cdot\kern-.33em\triangleright\!\!\!<}}
\def\lbiprod{{>\!\!\!\triangleleft\kern-.33em\cdot\, }}

\def\tens{\mathop{\otimes}}

\def\la{{\triangleright}}\def\ra{{\triangleleft}}

\def\isom{{\cong}}

\def\Hom{{\rm Hom}}
\def\Nat{{\rm Nat}}

\def\ev{{\rm ev}}
\def\coev{{\rm coev}}

\def\image{{\rm image}\, }

\def\id{{\rm id}}
\def\Deltaop{{\Delta^{\rm op}}}
\def\Lin{{\rm Lin}}

\def\<{\langle}
\def\>{\rangle}

\def\equad{\kern -1.7em}

\def\nquad{{\!\!\!\!\!\!}}
\def\nqquad{\nquad\nquad}
\def\eqn#1#2{\begin{equation}#2\label{#1}\end{equation}}

\def\haj#1{{\mathaccent20 {#1}}}
\def\Vhaj{{V\haj{\ }}}
\def\Whaj{{W\haj{\ }}}

\def\o{{}_{\scriptscriptstyle(1)}}
\def\t{{}_{\scriptscriptstyle(2)}}
\def\th{{}_{\scriptscriptstyle(3)}}

\def\bo{{}^{\bar{\scriptscriptstyle(1)}}}
\def\bt{{}^{\bar{\scriptscriptstyle(2)}}}
\def\Ro{{\CR^{\scriptscriptstyle(1)}}}
\def\Rt{{\CR^{\scriptscriptstyle(2)}}}

\def\und#1{{\underline {#1}}}

\def\uo{{{}^{\scriptscriptstyle(1)}}}
\def\ut{{{}^{\scriptscriptstyle(2)}}}

\def\Bo{{{}_{\und{\scriptscriptstyle(1)}}}}
\def\Bt{{{}_{\und{\scriptscriptstyle(2)}}}}

\def\text#1{\mbox{\rm #1}}
\def\note#1{}

\def\blacksquare{{\lform}}
\def\frac#1#2{{{#1\over#2}}}

\def\proof{\goodbreak\noindent{\bf Proof\quad}}

\def\endproof{{\ $\lform$}\bigskip }

\def\align#1{\begin{eqnarray*}#1\end{eqnarray*}}

\def\emi#1{{\em #1\index{#1}}}

\def\vect{{\bf t}}\def\vecv{{\bf v}}
\def\vecu{{\bf u}}\def\vecx{{\bf x}}

\def\span{{\rm span}}

\def\ant{{S}}
\documentstyle[12pt,epsf]{article}
\newtheorem{lemma}{Lemma}[section]
\newtheorem{propos}[lemma]{Proposition}
\newtheorem{example}[lemma]{Example}
\newtheorem{theorem}[lemma]{Theorem}

\newtheorem{corol}[lemma]{Corollary}
\newtheorem{defin}[lemma]{Definition}

\begin{document}
July 1993 \hskip 1.3in {\em Published in} Advances in Hopf Algebras,
Marcel Dekker

{\ }\hskip 2in  Lec. Notes Pure and Applied Maths {\bf 158} (1994) 55-105.
\vskip 1.3in
\baselineskip 11pt
\begin{center}{\Large ALGEBRAS AND HOPF ALGEBRAS\\ {\  }\\
IN BRAIDED CATEGORIES}\footnote{1991
Mathematics Subject Classification 18D10, 18D35, 16W30, 57M25, 81R50, 17B37
$\qquad\qquad\qquad$
This paper is in final form and no version of it will be submitted for
publication elsewhere} \\
{\ }\\ {\ }\\
{\ }\\ {\small SHAHN MAJID}\footnote{SERC Fellow and Fellow of Pembroke
College, Cambridge}\\ {\ }\\
{\it Department of Applied Mathematics \& Theoretical Physics}\\ {\it
University of Cambridge, Cambridge CB3 9EW, U.K.}
\end{center}

\begin{quote}
\vskip 15pt
\centerline{\small ABSTRACT} \small
This is an introduction for algebraists to the theory of algebras and Hopf
algebras in braided
categories. Such objects generalise super-algebras and super-Hopf algebras, as
well as colour-Lie algebras. Basic facts about braided categories $\CC$ are
recalled, the
modules and comodules of Hopf algebras {\em in} such categories are studied,
the notion of `braided-commutative' or `braided-cocommutative' Hopf algebras
(braided groups) is
reviewed and a fully diagrammatic proof of the reconstruction theorem for a
braided group $\aut(\CC)$ is given. The theory has important implications for
the theory of quasitriangular Hopf algebras (quantum groups). It also includes
important examples such as the degenerate Sklyanin algebra and the quantum
plane.
\end{quote}
 \baselineskip 15pt

One of the main motivations of the theory of Hopf algebras is that they
provide
a generalization of groups. Hopf
algebras of functions on groups provide examples of commutative Hopf algebras,
but it turns out that many group-theoretical constructions work just as well
when the Hopf algebra is allowed to be non-commutative. This is the philosophy
associated to some kind of non-commutative (or so-called quantum) algebraic
geometry. In a Hopf algebra context one can say the same thing in a dual way:
group algebras and enveloping algebras are cocommutative but many
constructions
are not tied to this. This point of view has been highly successful in recent
years, especially in
regard to the quasitriangular Hopf algebras of Drinfeld\cite{Dri}. These are
non-cocommutative but the non-cocommutativity is controlled by a
quasitriangular structure $\CR$. Such objects are commonly called quantum
groups. Coming out of physics, notably associated to solutions of the Quantum
Yang-Baxter Equations (QYBE) is a rich supply of
quantum groups.

Here we want to describe some kind of rival or variant of these quantum
groups,
which we call braided
groups\cite{Ma:bra}--\cite{Ma:skl}. These
are motivated by an earlier revolution
that was very popular some decades ago in mathematics and physics, namely the
theory of super or $\Z_2$-graded algebras
and Hopf algebras. Rather than make the algebras non-commutative etc one makes
the notion of tensor product $\tens$ non-commutative. The algebras remain
commutative with respect to this new tensor product (they are
super-commutative). Under this point of view one has super-groups,
super-manifolds and super-differential geometry. In many ways this line of
development was somewhat easier than the notion of quantum geometry because it
is conceptually easier to make an entire shift of category from vector spaces
to super-vector spaces. One can study Hopf algebras in such categories also
(super-quantum groups).

In this second line of development an obvious (and easy) step was to
generalise
such constructions to the case
of symmetric tensor categories\cite{Mac:cat}. These have a tensor product
$\tens$ and a collection of isomorphisms $\Psi$ generalizing the transposition
or super-transposition map but retaining its general properties. In
particular,
one keeps $\Psi^2=\id$ so that these generalized transpositions still generate
a representation of the symmetric group. Since only such general properties
are
used in most algebraic constructions, such as Hopf algebras and Lie algebras,
these notions immediately (and obviously) generalise to this setting. See for
example Gurevich \cite{Gur:yan}, Pareigis \cite{Par:non}, Scheunert
\cite{Sch:gen}
and numerous other authors. On the other hand, the theory is {\em not}
fundamentally
different from the super-case.

Rather more interesting is the further generalization to relax the condition
that $\Psi^2=\id$. Now $\Psi$ and $\Psi^{-1}$ must be distinguished and are
more conveniently represented by braid-crossings rather than by permutations.
They generate an action of the Artin braid group on tensor products. Such
quasitensor or braided-tensor categories have been formally introduced into
category theory in \cite{JoyStr:bra} and also arise in the representation
theory of quantum groups. The study of algebras and Hopf algebras in such
categories is rather more non-trivial than in the symmetric case. It is this
theory that we wish to describe here. It has been introduced by the author
under the heading `braided groups' as mentioned above. Introduced were the
relevant notions (not all of them obvious), the basic lemmas (such as a
braided-tensor product analogous to the super-tensor product of
super-algebras)
and a construction leading to a rich supply of examples.

On the mathematical level this project of `braiding' all of mathematics is, I
believe, a deep one (provided one goes from the symmetric to the truly braided
case). Much of mathematics consists of manipulating symbols, making
transpositions etc. The situation appears to be that in many constructions the
role of permutation group can (with care) be played equally well by the braid
group. Not only the algebras and braided groups to be described here, but also
braided differential calculus, braided-binomial theorems and
braided-exponentials are known\cite{Ma:fre} as well as braided-Lie
algebras\cite{Ma:lie}.  Much more can be expected. Ultimately we would like
some kind of braided geometry comparable to the high-level of development in
the super case.

Apart from this long-term philosophical motivation, one can ask what are the
more immediate applications of this kind of braided geometry? I would like to
mention five of them.

\begin{enumerate}
\item Many algebras of interest in physics such as the degenerate Sklyanin
algebra, quantum planes and exchange algebras are not naturally quantum groups
but turn out to be braided ones\cite{Ma:skl}\cite{Ma:poi}. There are
braided-matrices $B(R)$ and braided-vectors $V(R')$ associated to
$R$-matrices.

\item The category of Hopf algebras is not closed under quotients in a good
sense. For example, if $H\subset H_1$ is covered by a Hopf algebra projection
then $H_1\isom B\cocross H$ where $B$ is a braided-Hopf algebra. This is
the right setting for Radford's theorem as we have discovered and
explained in detail in \cite{Ma:skl}.

\item Braided groups are best handled by means of braid diagrams in which
algebraic operations `flow' along strings. This means deep connections with
knot theory and is also useful even for ordinary Hopf algebras. For example,
you can dualise theorems geometrically by turning the diagram-proof
up-side-down and flip conventions by viewing in  a mirror.

\item A useful tool in the theory of quasitriangular Hopf algebras (quantum
groups) via a process of {\em transmutation}. By encoding their
non-cocommutativity as braiding in a braided category they appear
`cocommutative'. Likewise, dual quasitriangular Hopf algebras are rendered
`commutative' by this process \cite{Ma:tra}\cite{Ma:bg}.

\item In particular, properties of the quantum groups $\CO_q(G)$ and $U_q(g)$
are
most easily understood in terms of their braided versions $B_q(G)$ and
$BU_q(g)$. This includes an Ad-invariant `Lie
algebra-like' subspace $\CL\subset U_q(g)$ and an isomorphism $B_q(G)\isom
U_q(g)$ \cite{Ma:exa}\cite{Ma:skl}.
\end{enumerate}

An outline of the paper is the following. In Section~1 we recall the basic
notions of braided tensor categories and how to work in them, and some
examples. We recall basic facts about quasitriangular and
dual quasitriangular Hopf algebras and the braided categories they generate.
In
Section~2 we do diagrammatic Hopf-algebra theory in this setting. In Section~3
we give a new diagrammatic proof of our generalised Tannaka-Krein-type
reconstruction theorem. In Section~4 we explain the results about ordinary
quantum groups obtained from this braided theory. In Section~5 we end with
basic examples of braided matrices etc associated to an $R$-matrix. Although
subsequently of interest in physics, the braided matrices arose quite
literally
from the Tannaka-Krein theorem mentioned above. This is an example of pure
mathematics feeding back into physics rather than the other way around
(for a change).

Our work on braided groups (or Hopf algebras in braided categories) was
presented to the Hopf algebra community at the Euler Institute in Leningrad,
October 1990 and
at the Biannual Meeting of the
American Maths Society in San Francisco, January 1991 and published in
\cite{Ma:eul}\cite{Ma:bg}.
The result presented at these meetings was the introduction of Hopf algebras
living {\em in} the
braided category of comodules of a dual quasitriangular Hopf
algebra. The connection between crossed modules (also
called Drinfeld-Yetter categories) and the quantum double as well as the
connection with Radford's theorem were introduced in \cite{Ma:dou} in early
1990. The braided interpretation
of Radford's theorem was introduced in detail in \cite{Ma:skl} and circulated
at the start of 1992.
Dual quasitriangular (or coquasitriangular) Hopf algebras
themselves were developed in connection with Tannaka-Krein ideas in \cite[Sec.
4]{Ma:pro}\cite{Ma:eul}\cite[Appendix]{Ma:bg} (and earlier in
other equivalent forms). A related Tannaka-Krein theorem in the
quasi-associative dual quasitriangular setting was obtained in
\cite{Ma:tan} at the Amherst conference and circulated in final form in
the Fall of 1990.

It is a pleasure to see that some of these ideas have subsequently proven of
interest in Hopf algebra circles (directly or indirectly). I would also like
to
mention some constructions of Lyubashenko\cite{Lyu:tan}\cite{Lyu:mod}
relating to our joint
work\cite{LyuMa:bra}. Also in joint work with Gurevich\cite{GurMa:bra} the
transmutation construction is related to Drinfeld's process of
twisting\cite{Dri:quas}.
Several other papers can be mentioned here. On the whole I have resisted the
temptation to give a full survey of all results obtained so far. Instead, the
aim here is a more pedagogical exposition of the more elementary results, with
proofs.

Throughout this paper we assume familiarity with usual techniques of Hopf
algebras such as in the book of Sweedler\cite{Swe:hop}. In this sense the
style
(and also the motivation) is somewhat different from our
braided-groups review article for physicists\cite{Ma:introp}. We work over a
field $k$. With more care one can work here with a ring just as well. When
working
with matrix or tensor components we will use the convention of summing over
repeated
indices. Some of the elementary quantum groups material should appear in more
detail in
my forthcoming book.

\section{Braided Categories}

Here we develop the braided categories within which we intend to work, namely
those coming
from (co)modules of quantum groups. In fact, the theory in Sections 2,3 is not
tied to quantum groups and works in any braided category. The material in the
present section is perfectly standard by now.

\subsection{Definition and General Constructions}

Symmetric monoidal (=tensor) categories have been known for some time and we
refer to \cite{Mac:cat} for
details. The model is the category of $k$-modules. The notion of braided
monoidal (=braided tensor=quasitensor) category is a small generalization if
this.

Briefly, a {\em monoidal category} means $(\CC,\tens,\Phi_{V,W,Z},\und 1,l,r)$
where $\CC$ is a category with objects $V,W,Z$ etc, $\tens:\CC\times \CC\to
\CC$
is a functor and $\Phi$ is a natural transformation between the two functors
$\ \tens(\ \tens\ )$ and $(\ \tens\ )\tens\ $ from $\CC\times\CC\times\CC\to
\CC$. This means a functorial collection of isomorphisms $\Phi_{V,W,Z}:V\tens
(W\tens Z)\to (V\tens W)\tens Z$.
These are in addition required to obey the `pentagon' coherence condition of
Mac Lane. This expresses equality of two ways to go via $\Phi$ from $U\tens
(V\tens(W\tens Z))\to ((U\tens V)\tens W)\tens Z$. Once this is assumed Mac
Lane's theorem ensures that all other re-bracketing operations are consistent.
In practice this means we can forget $\Phi$ and brackets entirely. We also
assume a unit object $\und1$ for the tensor product and associated functorial
isomorphisms $l_V:V\to \und 1\tens V,r_V:V\to V\tens \und 1$ for all objects
$V$, which we likewise suppress.

A monoidal category $\CC$ is {\em rigid} (=has left duals) if for each object
$V$, there is an object $V^*$ and morphisms
$\ev_V:V^*\tens V\to \und 1$, $\coev_V:\und 1
\to V\tens V^*$ such that
\eqn{lduala}{ V{\buildrel \coev
\over\to}(V\tens V^*)\tens V
\isom V\tens (V^*\tens V){\buildrel \ev\over \to}V}
\eqn{ldualb}{ V^*{\buildrel\coev
\over\to}V^*\tens (V\tens V^*)
\isom (V^*\tens V)\tens V^*{\buildrel \ev\over \to}
V^*}
compose to $\id_V$ and $\id_{V^*}$ respectively. A single object has a left
dual if $V^*,\ev_V,\coev_V$ exist. The model is that of a finite-dimensional
vector space (or finitely generated projective module when $k$ is a ring).

Finally, the monoidal category is {\em braided} if it has a {\em
quasisymmetry} or `braiding' $\Psi$ given as a natural transformation between
the two functors $\tens$ and $\tens^{\rm op}$ (with opposite product) from
$\CC\times\CC\to\CC$. This is a collection of functorial isomorphisms
$\Psi_{V,W}:V\tens W\to W\tens V$ obeying two `hexagon'
coherence identities. In our suppressed notation these are
\eqn{Psi-hex}{\Psi_{V\tens W,Z}=\Psi_{V,Z}\circ\Psi_{W,Z},\qquad
\Psi_{V,W\tens
Z}=\Psi_{V,Z}\circ\Psi_{V,W}}
while identities such as
\eqn{Psi-1}{\Psi_{V,\und 1}=\id_V=\Psi_{\und 1,V}}
can be deduced. If $\Psi^2=\id$ then one of the
hexagons is superfluous and we have an ordinary symmetric monoidal category.

Let us recall that the functoriality of maps such as those above means that
they commute
in a certain sense with morphisms in the category. For example, functoriality
of $\Psi$ means
\eqn{Psi-funct}{\Psi_{Z,W}(\phi\tens\id)=(\id\tens\phi)\Psi_{V,W}\
\forall\phi\matrix{\scriptstyle V\cr\downarrow\cr \scriptstyle Z},\qquad
\Psi_{V,Z}(\id\tens\phi)=(\phi\tens\id)\Psi_{V,W}\ \forall
\phi\matrix{\scriptstyle W\cr\downarrow\cr \scriptstyle Z}.}

These conditions (\ref{Psi-hex})-(\ref{Psi-funct}) are just the obvious
properties that we take for granted when transposing ordinary vector spaces or
super-vector spaces. In these cases $\Psi$ is the twist map $\Psi_{V,W}(v\tens
w)=w\tens v$ or the supertwist
\eqn{Psi-sup}{\Psi_{V,W}(v\tens w)=(-1)^{|v| |w|} w\tens v}
on homogeneous elements of degree $|v|,|w|$. The form of $\Psi$ in these
familiar cases does not depend directly on the spaces $V,W$ so we often forget
this. But in principle there is a different map $\Psi_{V,W}$ for each $V,W$
and
they all connect together as explained.

In particular, note that for any two $V,W$ we have two morphisms
$\Psi_{V,W},\Psi_{W,V}^{-1}:V\tens W\to W\tens V$ and in the truly braided
case
these can be distinct. A convenient notation in this case is to write them not
as permutations but as braid crossings. Thus we write morphisms pointing
downwards (say) and instead of a usual arrow, we use the shorthand
\eqn{Psi-bra}{\epsfbox{psi.eps}.}
In this notation the hexagons (\ref{Psi-hex}) appear as
\eqn{Psi-hex-bra}{\epsfbox{hexagons.eps}.}
The doubled lines refer to the composite objects $V\tens W$ and $W\tens Z$ in
a
convenient extension of the notation. The coherence theorem for braided
categories
can be stated very simply in this
notation: if two series of morphisms
built from $\Psi,\Phi$ correspond to the same braid then they compose to the
same morphism. The proof is just the same as Mac Lane's proof in the symmetric
case with the action of the symmetric group replaced by that of the Artin
braid
group.

This notation is a powerful one. We can augment it further by writing any
other
morphisms as nodes on a string connecting the inputs down to the outputs.
Functoriality (\ref{Psi-funct}) then says that a morphism $\phi:V\to Z$ say
can
be pulled through braid crossings,
\eqn{Psi-funct-bra}{\epsfbox{functorial.eps}}
Similarly for $\Psi^{-1}$ with inverse braid crossings. An easy lemma using
this
notation is that for any braided category $\CC$
there is another mirror-reversed braided monoidal category $\bar{\CC}$ with
the
same monoidal structure but with braiding
\eqn{psi-barC}{\bar{\Psi}_{V,W}=\Psi_{W,V}^{-1}}
in place of $\Psi_{V,W}$, i.e with the interpretation of braid crossings and
inverse braid crossings interchanged.

Finally, because of (\ref{Psi-1}) we can suppress the unit object entirely so
the evaluation and co-evaluation appear simply as $\ev=\cup$ and $\coev=\cap$.
Then (\ref{lduala})-(\ref{ldualb}) appear as
\eqn{lduals-bra}{\epsfxsize=5.7in \epsfbox{lduals.eps}.}
There is a similar notion of right duals $\Vhaj$ and $\bar\ev_V, \bar\coev_V$
for
which the mirror-reflected
double-bend here can be likewise straightened.

\begin{example} Let $R\in M_n(k)\tens M_n(k)$ be invertible and obey the QYBE
\[ R_{12}R_{13}R_{23}=R_{23}R_{13}R_{12}\]
then the monoidal category $\CC(V,R)$ generated by tensor products of $V=\C^n$
is braided.
\end{example}
\proof This is an elementary exercise (and extremely well-known).  The
notation
is $R_{12}=R\tens\id$ and $R_{23}=\id\tens R$ in $M_n^{\tens 3}$. The braiding
on basis vectors $\{e_i\}$ is
\eqn{Psi-vec}{\Psi(e_i\tens e_j)=e_b\tens e_a R^a{}_i{}^b{}_j}
extended to tensor products according to (\ref{Psi-hex}). The morphisms in the
category are linear maps such that $\Psi$ is functorial with respect to them
in
the sense of (\ref{Psi-funct}). The associativity $\Phi$ is the usual one on
vector spaces. \endproof

If $R$ obeys further conditions then $\CC(V,V^*,R)$ generated by $V,V^*$ is
rigid.
One says that
such an $R$  is {\em dualizable}. For this there should exist among other
things a
`second-inverse'
\eqn{tildeR}{\widetilde R=((R^{t_2})^{-1})^{t_2}}
where $t_2$ is transposition in the second $M_n$ factor. This defines one of
the mixed terms in the braiding
\eqn{Psi-co-co}{\Psi_{V^*,V^*}(f^i\tens f^j)=R^i{}_a{}^j{}_b f^b\tens f^a}
\eqn{Psi-vec-co}{\Psi_{V,V^*}(e_i\tens f^j)=\widetilde{R}^a{}_i{}^j{}_b
f^b\tens e_a}
\eqn{Psi-co-vec}{\Psi_{V^*,V}(f^i\tens e_j)=e_a\tens f^b
R^{-1}{}^i{}_b{}^a{}_j}
where $V^*=\{f^i\}$ is a dual basis. The evaluation and coevaluation are given
by the usual morphisms
\eqn{ev-coev}{ \ev_V(f^i\tens e_j)=\delta^i{}_j,\quad \coev_V(1)=\sum_i
e_i\tens f^i.}
One needs also the second-inverse $\widetilde{R^{-1}}$ for $\Psi$ to be
invertible.
In this way one translates the various axioms into a linear space
setting. We see in particular that the QYBE are nothing other than the braid
relations in matrix form.

We turn now to some general categorical constructions. One construction in
\cite{Ma:rep}\cite{Ma:cat} is based on the idea that a pair of monoidal
categories $\CC\to\CV$ connected by a functor behaves in many ways like a
bialgebra with $\tens$ in $\CC$ something like the product. In some cases
this is actually true as we shall see in Section~3 in the form of a
Tannaka-Krein-type reconstruction theorem, but we can keep it in general as
motivation. Motivated by this we showed that for every pair $\CC\to\CV$ of
monoidal categories there is a dual one $\CC^\circ\to\CV$ where $\CC^\circ$ is
the {\em Pontryagin dual monoidal category}\cite{Ma:rep}. This generalised the
usual duality for Abelian groups and bialgebras to the setting of monoidal
categories. We also proved such things as a canonical functor
\eqn{pontfunc}{\CC\to {}^\circ(\CC^\circ).}
Of special interest to us now is the case $\CC\to\CC$ where the functor is the
identity one. So associated to every monoidal category $\CC$ is another
monoidal
category $\CC^\circ$ of
`representations' of $\tens$. This special case can also be denoted by
$\CC^\circ=\CZ(\CC)$ the `center' or `inner double' of $\CC$ for reasons that
we
shall explain shortly. This case was found independently by V.G. Drinfeld who
pointed out that it is braided.

\begin{propos}\cite{Ma:rep}\cite{Dri:pri} Let $\CC$ be a monoidal category.
There is a braided monoidal category $\CC^\circ=\CZ(\CC)$ defined as
follows. Objects are pairs $(V,\lambda_V)$ where $V$ is an object of $\CC$ and
$\lambda_V$ is a natural isomorphism in $\Nat(V\tens\id,\id\tens V)$ such that
\[ \lambda_{V,\und 1}=\id,\qquad
(\id\tens\lambda_{V,Z})(\lambda_{V,W}\tens\id)=\lambda_{V,W\tens Z}. \]
and morphisms are $\phi:V\to W$ such that the modules are intertwined in the
form
\[ (\id\tens\phi)\lambda_{V,Z}=\lambda_{W,Z}(\phi\tens\id),\qquad \forall Z\
{\rm in}\ \CC.\]
The monoidal product and braiding are
\[ (V,\lambda_V)\tens (W,\lambda_W)=(V\tens W,\lambda_{V\tens W}),\
\lambda_{V\tens W,Z}=(\lambda_{V,Z}\tens\id)(\id\tens\lambda_{W,Z})\]
\[ \Psi_{(V,\lambda_V),(W,\lambda_W)}=\lambda_{V,W}.\]
\end{propos}
\proof The monoidal structure was found in the author's paper \cite{Ma:rep}
where full proofs were also given. We refer to this for details. Its preprint
was circulated in
the Fall of 1989.
The braiding was pointed out by
Drinfeld\cite{Dri:pri} who had considered the construction from a very
different and
independent point of view to our duality one, namely in connection with the
double of a Hopf algebra as we shall explain below.  Another claim to the
construction is from the direction of tortile categories\cite{JoyStr:tor}.
See also \cite{Ma:cat} for further work from
the duality point of view.
\endproof

The `double' point of view for this construction is based on the following
example cf\cite{Dri:pri}.

\begin{example} Let $H$ be a bialgebra over $k$ and $\CC={}_H\CM$ the monoidal
category of $H$-modules. Then an object of $\CZ(\CC)$ is a vector space $V$
which
is both a left $H$-module and an invertible left $H$-comodule such that
\[  \sum h\o v\bo\tens h\t\la v\bt=\sum (h\o\la v)\bo h\t\tens (h\o\la v)\bt,
\qquad\forall h\in H,\ v\in V. \]
In this form $\CZ({}_H\CM)$ coincides with the category ${}_H^H\CM$ of
$H$-crossed modules\cite{Yet:rep} with an additional invertibility condition.
The
braiding is
\[ \Psi_{V,W}(v\tens w)=\sum v\bo\la w\tens v\bt. \]
The invertibility condition on the comodules ensures that $\Psi^{-1}$ exists,
and
is automatic if the bialgebra $H$ has a skew-antipode.
\end{example}
\proof The proof is standard from the point of view of Tannaka-Krein
reconstruction
methods (which we shall come to later). From $\CC$ we can reconstruct $H$
as the representing object for a certain functor. This establishes a bijection
$\Lin(V,H\tens V)\isom \Nat(V\tens \id,\id\tens V)$ under which $\lambda_V$
corresponds to a map $V\to H\tens V$. That $\lambda_V$ represents $\tens$
corresponds then to the comodule property of this map. That $\lambda_V$ is a
collection of morphisms corresponds to the stated compatibility condition
between this coaction and the action on $V$ as an object in $\CC$. To see this
in detail let $H_L$ denote $H$ as an object in $\CC$ under the left action.
Given $\lambda_V$ a natural transformation we define
\eqn{lamcoact}{ \sum v\bo\tens v\bt=\lambda_{V,H_L}(v\tens 1)}
and check
\align{(\id\tens\lambda_{V,H_L})(\lambda_{V,H_L}\tens\id)(v\tens 1\tens 1) &=&
\lambda_{V,H_L\tens H_L}(v\tens(1\tens 1))\\
&=&\lambda_{V,H_L\tens H_L}(v\tens \Delta(1))=(\Delta\tens
\id)\circ\lambda_{V,H_L}(v\tens 1)}
where the first equality is the fact that $\lambda_V$ `represents' $\tens$ and
the
last is that $\lambda_V$ is functorial under the morphism $\Delta:H_L\to
H_L\tens
H_L$. The left hand side is the map $V\to H\tens V$ in (\ref{lamcoact})
applied
twice so we see that this map is a left coaction. Moreover,
\align{\sum h\o v\bo\tens h\t\la v\bt &=& h\la\lambda_{V,H_L}(v\tens
1)=\lambda_{V,H_L}(h\la(v\tens 1))\\
&=&\sum \lambda_{V,H_L}(h\o\la v\tens R_{h\t}(1))=\sum \left(\lambda_{V,H_L}
(h\o\la v\tens
1)\right)(h\t\tens 1)}
where the first equality is the definition (\ref{lamcoact}) and the action of
$H$ on $H_L\tens V$. The second equality is that $\lambda_{V,H_L}$ is a
morphism in $\CC$. The final equality uses functoriality under the morphism
$R_{h\t}:H_L\to H_L$ given by right-multiplication to obtain the right hand
side
of the compatibility
condition. The converse directions are easier. Given a coaction $V\to H\tens
V$
define $\lambda_{V,W}(v\tens w)=\sum v\bo\la w\tens v\bt$. This also implies
at
once the braiding
$\Psi=\lambda$ as stated.

Finally we note that in Proposition~1.2 the definition assumes that the
$\lambda_V$
are invertible. If we were to relax this then we would have a monoidal
category
which
is just that of crossed modules as in \cite{Yet:rep}, but then $\Psi$ would
not
necessarily be invertible and hence would not be a true braiding. The
invertible
$\lambda_V$ correspond to left comodules which are {\em invertible} in the
following
sense: there exists a linear map $V\to V\tens H$ sending $v$ to $\sum
v^{[2]}\tens v^{[1]}$ say, such that
\eqn{inv-comod}{\sum v^{[2]}\bo v^{[1]}\tens v^{[2]}\bt=1\tens v=\sum
v\bt{}^{[1]}v\bo\tens v\bt{}^{[2]},\qquad \forall v\in V.}
One can see that if such an `inverse' exists, it is unique and a right
comodule.
Moreover, it is easy to see that the invertible comodules are closed under
tensor
products. They correspond to $\lambda_V^{-1}$ in a similar way to
(\ref{lamcoact}) and with $\lambda_{V,W}^{-1}(w\tens v)=\sum v^{[2]}\tens
v^{[1]}\la w$ for the converse direction.
In the finite-dimensional case they provide left duals $V^*$
with left coaction $\beta_{V^*}(f)(v)=\sum v^{[1]}f(v^{[2]})$. If the
bialgebra
$H$ has a skew-antipode then
every left comodule is invertible by composing with the skew-antipode. So in
this
case the condition becomes empty.

{}From the categorical point of view in Proposition~1.2, if $\CC$ has right
duals then every $\lambda_{V,W}$ is invertible, cf\cite{Ma:rep}. The inverse
is
the right-adjoint of $\lambda_{V,\Whaj}$, namely
$\lambda_{V,W}^{-1}=(\bar\ev_W\tens\id)\circ\lambda_{V,\Whaj}\circ(\id\tens
\bar\coev_W)$. When  $\CC={}_H\CM$ then the finite-dimensional left modules
have right duals if the bialgebra $H$ has a skew-antipode, so in this case the
invertibility of $\lambda_V$ is automatic. On the other hand, we do not need
to
make these suppositions here.

 This completes our computation of $\CZ({}_H\CM)$. Apart from the
invertibility
restriction we see that it consists of compatible module-comodule structures
as
stated. \endproof

Note that the notion of a crossed module is an immediate generalisation of the
notion of a crossed G-module\cite{Whi:com} with $H=kG$, the group algebra of a
finite group $G$. In this case the category of crossed $G$-modules is
well-known to be braided\cite{FreYet:bra}. Moreover,  because the objects can
be identified with underlying vector spaces, we know by the Tannaka-Krein
reconstruction theorem\cite{Sav:cat} that
there must exist a bialgebra $coD(H)$ such that our braided-category is
equivalent to that of right $coD(H)$-comodules
\eqn{bialgdouble}{\CM_{f.d.}^{coD(H)}={}_H^H\CM_{f.d.}}
Here we take the modules to be finite-dimensional as a sufficient (but not
necessary) condition for a Tannaka-Krein reconstruction theorem to apply and
the co-double $coD(H)$ to exist. In the nicest case the category is also
${}_{D(H)}\CM_{f.d.}$ for some $D(H)$. This is an abstract definition of
Drinfeld's quantum double and works for a bialgebra.

If it happens that $H$ is a Hopf algebra with invertible antipode then one can
see from the above that ${}_H^H\CM_{f.d.}$ is rigid and so $coD(H)$ and $D(H)$
will be Hopf algebras. The categorical reason is that ${}_H\CM_{f.d.}$ is
rigid
and this duality extends to $\CZ(\CC)$ with the dual of $\lambda_V$ defined by
the left-adjoint of $\lambda_{V,W}^{-1}$, namely
$\lambda_{V^*,W}=(\ev_V\tens\id)\circ\lambda_{V,W}^{-1}\circ(\id\tens
\coev_V)$. We will
study details about categories of modules and comodules and the reconstruction
theorems later in this section and in Section~3. The point is that these
categorical methods are  very powerful.

\begin{propos}\cite{Ma:phy}cf\cite{Dri} If $H$ is a finite-dimensional Hopf
algebra then $D(H)$ (the quantum double Hopf algebra of $H$) is built on
$H^*\tens H$ as a
coalgebra with the product
\[ (a\tens h)(b\tens g)=\sum b\t a\tens h\t g<Sh\o,b\o><h\th,b\th>,\ h,g\in
H,\, a,b\in H^*\]
where $<\ ,\ >$ denotes evaluation.
\end{propos}
\proof The quantum double $D(H)$ was introduced by Drinfeld\cite{Dri} as a
system of generators and relations built from the structure constants of $H$.
The formula stated on $H^*\tens H$ is easily obtained from this as done in
\cite{Ma:phy}. We have used here the conventions introduced in \cite{Ma:dou}
that avoid the use
of the  inverse of the antipode. Also in \cite{Ma:dou} we showed that the
modules of the double were precisely
the crossed modules category as required. To see this simply note that $H$ and
$H^{*\rm op}$ are sub-Hopf algebras and hence a left $D(H)$-module is a left
$H$-module and a suitably-compatible right $H^{*\rm}$-module. The latter is
equally well a left $H$-comodule compatible as in Example~1.3. See
\cite{Ma:dou} for details. \endproof

In \cite{Ma:phy} we introduced a further characterization of the quantum
double
as a member of a class of {\em double cross product} Hopf algebras $H_1\bowtie
H_2$ (in which $H_i$ are mutually acting on each other). Thus, $D(H)=H^{*\rm
op}\bowtie
H$ where the actions are mutual coadjoint actions. In this form it is
clear that the role of $H^*$ can be played by $H^\circ$ in the infinite
dimensional Hopf algebra case. We will not need this further here.

\subsection{Quasitriangular Hopf Algebras}

We have already described one source of braided categories, namely as modules
of the double $D(H)$ (or comodules of the codouble) of a bialgebra.
Abstracting
from this one has the notion, due to Drinfeld, of a quasitriangular Hopf
algebra. These are such that their category of modules is braided.

\begin{defin}\cite{Dri} A \emi{quasitriangular bialgebra or Hopf algebra} is a
pair $(H,\CR)$ where $H$ is a bialgebra or Hopf algebra and $\CR\in H\tens
H$ is
invertible and obeys
\eqn{qua1}{(\Delta\tens \id)\CR=\CR_{13}\CR_{23},\quad
(\id\tens\Delta)\CR=\CR_{13}\CR_{12}.}
\eqn{qua2}{\tau\circ\Delta h=\CR(\Delta h)\CR^{-1}, \ \forall h\in H.}
Here $\CR_{12}=\CR\tens 1$ and $\CR_{23}=1\tens\CR$ etc, and $\tau$ is the
usual twist map.
\end{defin}

Thus these Hopf algebras are like cocommutative enveloping algebras or group
algebras but are cocommutative now only up to an isomorphism implemented by
conjugation by an element $\CR$. Some elementary (but important) properties
are

\begin{lemma}\cite{Dri:alm} If $(H,\CR)$ is a quasitriangular bialgebra then
$\CR$ as an element of $H\tens H$ obeys
\eqn{quaeps}{(\eps\tens \id)\CR=(\id\tens \eps)\CR=1.}
\eqn{uniqybe}{\CR_{12}\CR_{13}\CR_{23}=\CR_{23}\CR_{13}\CR_{12}}
If $H$ is a Hopf algebra then one also has
\eqn{quaant}{(S\tens\id)\CR=\CR^{-1},\qquad (\id\tens
S)\CR^{-1}=\CR,\quad (S\tens S)\CR=\CR}
\eqn{uv}{\exists S^{-1},u,v,\quad S^2(h)=u h u^{-1},\quad S^{-2}(h)=v h
v^{-1}\quad \forall h\in H}
\end{lemma}
\proof For (\ref{quaeps}) apply $\eps$ to (\ref{qua1}), thus
$(\eps\tens\id\tens\id)(\Delta\tens\id)\CR=\CR_{23}=(\eps\tens\id\tens\id)
\CR_{13}\CR_{23}$ so that (since $\CR_{23}$ is invertible) we have
$(\eps\tens\id)\CR=1$. Similarly for the other side. For (\ref{uniqybe})
compute $(\id\tens\tau\circ\Delta)\CR$ in two ways: using the second of
axioms (\ref{qua1}) directly
or using axiom (\ref{qua2}), and then the second of (\ref{qua1}). For
(\ref{quaant}) consider $\sum \Ro\o S
\Ro\t\tens\Rt=1$ by the property of the antipode and equation (\ref{quaeps})
already proven, but
equals $\CR(S\tens\id)\CR$ by axiom (\ref{qua1}). Similarly for the other
side, hence $(S\tens\id)\CR=\CR^{-1}$. Similarly for $(\id\tens
S)\CR^{-1}=\CR$
once we appreciate that
$(\Delta\tens\id)(\CR^{-1})=(\CR_{13}\CR_{23})^{-1}
=\CR_{23}^{-1}\CR_{13}^{-1}$
etc, since $\Delta$ is an algebra homomorphism. For (\ref{uv}) the relevant
expressions are
\[u=\sum (S\CR\ut)\CR\uo,\quad u^{-1}=\sum \CR\ut S^2\CR\uo,\quad v=Su\]
which one can verify to have the right properties. In addition one can see
that
$\Delta u=(\CR_{21}\CR_{12})^{-1}(u\tens u)$ and similarly for $v$ so that
$uv^{-1}$ is group-like (and implements $S^4$). For details of the
computations
see \cite{Dri:alm} or reviews by the author. \endproof

Here (\ref{uniqybe}) is the reason that Physicists call $\CR$ the `universal
R-matrix' (compare Example~1.1). Indeed, in any finite-dimensional
representation the image of $\CR$ is such an $R$-matrix. There are well-known
examples such as $U_q(sl_2)$ and $U_q(g)$\cite{Dri}\cite{Jim:dif}. Here we
give
perhaps the simplest known quasitriangular Hopf algebras

\begin{example}\cite{Ma:any} Let $\Z_n=\Z/n\Z$ be the finite cyclic group of
order $n$ and $k\Z_n$ its group algebra with generator $g$. Let $q$ be a
primitive $n$-th root of unity. Then there is a quasitriangular Hopf algebra
$\Z_n'$ consisting of this group algebra and
\eqn{any-R}{\Delta g=g\tens g,\quad \eps g=1,\quad S g=g^{-1},\quad
\CR=n^{-1}\sum_{a,b=0}^{n-1} q^{-ab}g^a\tens g^b.}
\end{example}
\proof We assume that $k$ is of suitable characteristic. To verify the
non-trivial quasitriangular structure we use  that
$n^{-1}\sum_{b=0}^{n-1}q^{ab}=\delta_{a,0}$. Then $\CR_{13}\CR_{23}=n^{-2}\sum
q^{-(ab+cd)}g^a\tens g^c\tens g^{b+d}$ $=$ $n^{-2}\sum q^{-
b(a-c)}q^{-cb'}g^a\tens g^c\tens g^{b'}$ $=$ $n^{-1}\sum q^{- ab'}g^a\tens
g^a\tens g^{b'}$ $=$ $(\Delta\tens\id)\CR$ where $b'=b+d$ was a change of
variables. Similarly for the second of (\ref{qua1}). The remaining axiom
(\ref{qua2}) is automatic because the Hopf algebra is both commutative and
cocommutative. \endproof

\begin{example}\cite{BMO:mul}\cite{Ma:csta} Let $G$ be a finite Abelian group
and $k(G)$ its
function Hopf algebra. Then a quasitriangular structure on $k(G)$ means a
function $\CR\in H\tens H$ obeying
\[ \CR(gh,f)=\CR(g,f)\CR(h,f),\quad \CR(g,hf)=\CR(g,h)\CR(g,f),\quad
\CR(g,e)=1=\CR(e,g)\]
for all $g,h,f$ in $G$ and $e$ the identity element. I.e., a quasitriangular
structure on $k(G)$ means precisely a bicharacter of $G$.
\end{example}
\proof We identify $k(G)\tens k(G)$ with functions on $G\times G$, with
pointwise multiplication. Using the comultiplication given by multiplication
in
$G$ we have at once that (\ref{qua1}) corresponds to the first two displayed
equations. Axiom (\ref{qua2}) becomes $hg\CR(g,h)=\CR(g,h)gh$ and so is
automatic because the group is Abelian. Given these first two of the stated
conditions, the latter two hold iff $\CR$ is invertible. \endproof

The $\Z_n'$ example here also has an immediate generalization to the
group algebra $k G$ of a finite Abelian group equipped with a bicharacter on
$\hat G$. This just coincides with the last example applied to $k(\hat G)=k
G$.
Finally, we return to our basic construction,

\begin{example}\cite{Dri} Let $H$ be a finite-dimensional bialgebra. Then
$D(H)$ is quasitriangular. In the Hopf algebra case the
quasitriangular structure is $\CR=\sum_a (f^a\tens 1)\tens (1\tens e_a)$ where
$H=\{e_a\}$ is a basis and $\{f^a\}$ a dual basis.
\end{example}
\proof The result is due to Drinfeld. A direct proof in the abstract Hopf
algebra setting appeared in \cite{Ma:phy}. The easiest way to show that
$\CR$ is invertible is to verify in view of (\ref{quaant}) that $(S\tens
\id)(\CR)$ is the inverse. \endproof

\begin{theorem}e.g.\cite{Ma:qua} Let $(H,\CR)$ be a quasitriangular bialgebra.
Then the category ${}_H\CM$ of modules is braided. In the Hopf algebra case
the
finite-dimensional modules are rigid. The braiding and the action on duals are
\[ \Psi_{V,W}(v\tens w)=\sum \CR\ut\la w\tens\CR\uo\la v,\quad  h\la
f=f((Sh)\la(\ ))\]
with $\ev,\coev$ as in (\ref{ev-coev}).
\end{theorem}
\proof $h\la\Psi(v\tens w)=(\Delta h)\la \tau(\CR\la
(v\tens w))=\tau((\Deltaop h)\CR\la (v\tens w))=\tau(\CR(\Delta h)\la (v\tens
w))=\Psi(h\la(v\tens w))$ in virtue (\ref{qua2}). It is easy to
see that (\ref{qua1}) likewise just correspond to the
hexagons (\ref{Psi-hex}) or (\ref{Psi-hex-bra}). Functoriality is also easily
shown. For an early treatment of this topic see \cite[Sec. 7]{Ma:qua}. Note
that if $\CR_{21}=\CR^{-1}$ (the triangular rather than quasitriangular case)
we have $\Psi$ symmetric rather than braided. This was the case treated in
\cite{Dri} though
surely the general quasitriangular case was also known to some experts at the
time or shortly thereafter. \endproof

\begin{propos}\cite[Sec. 6]{Ma:exa} Let $H=\Z_2'$ denote the quantum group
in Example~1.7 with $n=2$. Then $\CC={}_{\Z_2'}\CM={\rm SuperVec}$ the
category
of super-vector spaces.
\end{propos}
\proof One can easily check that this $\Z_2'$ is indeed a quasitriangular (in
fact, triangular) Hopf algebra. Hence we have a (symmetric) tensor category of
representations. Writing $p={1-g\over 2}$ we have $p^2=p$ hence any
representation $V$ splits as $V_0\oplus V_1$ according to the
eigenvalue of $p$. We can also write $\CR=1-2p\tens p$ and hence from
Theorem~1.10 we compute $\Psi(v\tens w)=\tau(\CR\la (v\tens w))=(1-2p\tens
p)(w\tens v)=(1-2|v||w|)w\tens v=(-1)^{|v||w|}w\tens v$ as in (\ref{Psi-sup}).
\endproof

So this non-standard quasitriangular Hopf algebra $\Z_2'$ (non-standard
because
of its
non-trivial $\CR$) recovers the category of super-spaces with its correct
symmetry $\Psi$. In just the same way the category $\CC_n={}_{\Z_n'}\CM$
consists of vector spaces that split as
$V=\oplus_{a=0}^{n-1} V_a$ with the degree of an element defined by the action
$g\la v=q^{|v|}v$ where $q$ is a primitive $n$-th root of unity. From
Theorem~1.10 and (\ref{any-R}) we
find
\eqn{Psi-any}{\Psi_{V,W}(v\tens w)=q^{\vert v\vert \vert w\vert}w\tens v.}
Thus we call $\CC_n$ the category of {\em anyonic vector spaces} of fractional
statistics ${1\over n}$, because just such a braiding is encountered in
anyonic
physics. For $n>2$ the category is
strictly braided in the sense that $\Psi\ne\Psi^{-1}$. There are natural
anyonic traces
and anyonic dimensions generalizing the super-case\cite{Ma:any}
\eqn{any-dim}{\und\dim (V)=\sum_{a=0}^{n-1} q^{-{a^2}}\dim
V_a,\qquad \und{\rm Tr}(f)=\sum_{a=0}^{n-1}q^{-{a^2}} \trace
f|_{V_a}.}
We see that this anyonic category is generated by the quantum group $\Z_n'$.

Obviously we can take this idea for generalising super-symmetry to the further
case of Example~1.8. In this case a $k(G)$-module just means a $G$-graded
space
where $f\la v=f(|v|)v$ on homogeneous elements of degree $|v|\in G$.
This is well-known to Hopf algebraists for some time: the new ingredient is
that
a bicharacter gives our $G$-graded spaces a natural braided-transposition
$\Psi$. We have given plenty of other less obvious examples of braided
categories generated in this way from quasitriangular Hopf algebras \cite[Sec.
6]{Ma:exa}\cite{Ma:tra}. Our idea in this work is not to use Hopf algebras in
connection with deformations (the usual setting) but rather as the `generator'
of a category within which we shall later make algebraic constructions. This
is
how quantum groups are naturally used to generalise supersymmetry. In this
context they are typically discrete.

\subsection{Dual Quasitriangular Structures}

In this section we describe the dual results to those above. If a
quasitriangular Hopf algebra is almost cocommutative up to conjugation then
its
dual Hopf algebra should be almost commutative up to `conjugation' in the
convolution algebra. The relevant axioms are obtained by dualizing in the
standard way by writing out the axioms of a quasitriangular Hopf algebra as
diagrams and then reversing all the arrows (and a left-right reversal).
Obviously it is the axioms that are being dualised and not any specific Hopf
algebra. This is important because in the infinite-dimensional case the dual
axioms are weaker. This is a rigorous way to work with the standard quantum
groups over a
field as appreciated in \cite{Ma:pro} among other places.

We will always denote our dual quasitriangular bialgebras and Hopf algebras by
$A$ (to avoid confusion). These are equipped now with a map $\CR:A\tens A\to
k$
which should be invertible in $\Hom(A\tens A,k)$ in the sense that there
exists
a map $\CR^{-1}:A\tens A\to k$ such that
\[\sum \CR^{-1}(a\o\tens b\o)\CR(a\t\tens b\t)=\eps(a)\eps(b)=\sum
\CR(a\o\tens
b\o)\CR^{-1}(a\t\tens
b\t).\]
Keeping such considerations in mind, it is easy to dualize the remainder of
Drinfeld's axioms to obtain the following definition.

\begin{defin}
 A dual quasitriangular bialgebra or Hopf algebra $(A,\CR)$ is a bialgebra
or
Hopf algebra $A$
and a convolution-invertible map $\CR:A\tens A\to k$ such that
\eqn{dqua1}{\CR(ab\tens c)=\sum
\CR(a\tens c\o)\CR(b\tens c\t),\quad \CR(a\tens bc)=\sum
\CR(a\o\tens c)\CR(a\t\tens b)}
\eqn{dqua2}{ \sum b\o a\o \CR(a\t\tens b\t)=\sum \CR(a\o\tens b\o) a\t b\t}
for all $a,b,c\in A$.
\end{defin}

This looks a little unfamiliar but is in fact obtained by replacing the
multiplication in Definition~1.5 by the convolution product and the
comultiplication by the multiplication in $A$. Axiom (\ref{dqua2}) is the dual
of (\ref{qua2}) and says, as promised, that $A$ is almost commutative -- up to
$\CR$. Axioms (\ref{dqua1}) are the dual of (\ref{qua1}) and say
that $\CR$ is a `bialgebra bicharacter'. They should be compared with
Example~1.14  below. We also have analogues of the various results in
Section~1.2. Again, the new language is perhaps unfamiliar so we give some of
the proofs in this dual form in detail.

\begin{lemma}\cite{Ma:bg} If $(A,\CR)$ is a dual quasitriangular bialgebra
then
\eqn{dquaeps}{ \CR(a\tens 1)=\eps(a)=\CR(1\tens a).}
\eqn{dqybe}{\sum \CR(a\o\tens b\o)\CR(a\t\tens c\o)\CR(b\t\tens c\t)=\sum
\CR(b\o\tens c\o)\CR(a\o\tens c\t)\CR(a\t\tens b\t)}
for all $a,b,c$ in $A$. If $A$ is a Hopf algebra then in addition,
\eqn{dquaant}{ \CR(Sa\tens b)=\CR^{-1}(a\tens b),\qquad \CR^{-1}(a\tens
Sb)=\CR(a\tens b),\quad \CR(Sa\tens Sb)=\CR(a\tens b)}
\eqn{duv}{\exists S^{-1},\quad v(a)=\sum \CR(a\o\tens Sa\t),\qquad  \sum a\o
v(a\t)=\sum v(a\o) S^2 a\t}
\end{lemma}
\proof Using (\ref{dqua1}) we have
\align{&&\nqquad\CR(a\tens 1)=\sum (\CR^{-1}(a\o\tens 1) \CR(a\t\tens 1))
\CR(a\th\tens 1)\\
&&=\sum \CR^{-1}(a\o\tens 1) (\CR(a\t\tens 1) \CR(a\th\tens 1))=\sum
\CR^{-1}(a\o\tens 1)
\CR(a\t\tens
1.1)=\eps(a)}
as in (\ref{dquaeps}). Likewise on the other side.  Also, if $\CR^{-1}$ exists
it is unique. Hence for $A$ a Hopf algebra it is given by $\CR^{-1}(a\tens
b)=\CR(Sa\tens b)$ (use axioms (\ref{dqua1})). In this case $a\tens
b\mapsto \CR(Sa\tens Sb)$ is convolution inverse to $\CR^{-1}$ because $
\sum \CR(Sa\o\tens Sb\o)$ $\CR(Sa\t\tens b\t)= \sum \CR(Sa\tens
(Sb\o)b\t)=\CR(Sa\tens
1)\eps(b)=\eps(a)\eps(b)$ etc. Hence $\CR(Sa\tens Sb)=\CR(a\tens b)$, proving
the other side and the third part of (\ref{dquaant}).
For (\ref{dqybe}) we apply the second of (\ref{dqua1}), (\ref{dqua2}) and the
second of (\ref{dqua1}) again,
\align{&&\nqquad\sum \left(\CR(a\o\tens b\o)\CR(a\t\tens
c\o)\right)\CR(b\t\tens
c\t)=\sum \CR(a\tens c\o b\o)\CR(b\t\tens c\t)\\
&&=\sum \CR(b\o\tens c\o)\CR(a\tens b\t c\t)=\sum \CR(b\o\tens
c\o)\left(\CR(a\o\tens
c\t)\CR(a\t\tens b\t)\right).}
For (\ref{duv}) one defines $v:A\to k$ as shown (and similarly a map $u:A\to
k$)
and checks the relevant facts analogous to Lemma~1.6. This is done in complete
detail in \cite[Appendix]{Ma:bg} to which we refer the reader. \endproof

Note that if $\CR$ is a linear map obeying (\ref{dquaeps}) and
(\ref{dqua1}) and if $A$ is a Hopf algebra, we can use
(\ref{dquaant}) as a definition of $\CR^{-1}$. Some authors in defining
similar
notions have made (\ref{dquaeps}) an axiom in the
bialgebra case (this is the case in \cite{Hay:det} and in first versions of
some other
works). For this and other reasons we stick to our original terminology
from\cite{Ma:pro}\cite{Ma:tan}\cite{Ma:eul}\cite{Ma:bg} with axioms and
properties as above.

\begin{example} Let $G$ be an Abelian group and $kG$ its group algebra.
This is dual quasitriangular iff there is a function $\CR:G\times G\to k$
obeying the bicharacter conditions in Example~1.8.
\end{example}
\proof In the group algebra we can work with group-like elements (these form a
basis). On such elements the axioms in (\ref{dqua1})-(\ref{dqua2}) simplify:
simply drop the ${}\o,\t$ suffixes! This immediately reduces to the
bimultiplicativity while invertibility corresponds once again (given this) to
(\ref{dquaeps}). \endproof

\goodbreak

This example is clearly identical in content to Example~1.8. Because $kG$ is
dual to $k(G)$ by evaluation, it is obvious that a dual quasitriangular
structure on $kG$ is just the same thing as a quasitriangular structure on
$k(G)$, namely as we see, a bicharacter. The result is however, more
transparent from the dual quasitriangular point of view and slightly more
general.
This example is the
reason that we called $\CR$ obeying (\ref{dqua1}) a bialgebra
bicharacter\cite{Ma:pro}. For a
concrete example, one can take $G=\Z_n$ and $\CR(a,b)=q^{ab}$ with $q$ a
primitive $n$-th root of unity to give a non-standard dual quasitriangular
structure on $k\Z_n$. It is just Example~1.7 after a $\Z_n$-Fourier transform.

Now let $R$ be an invertible matrix solution of the QYBE as in Example~1.1.
There is a by-now standard bialgebra $A(R)$\cite{Dri}\cite{FRT:lie} defined by
generators $1$ and $\vect=\{t^i{}_j\}$ (regarded as an $n\times n$ matrix) and
the relations, comultiplication and counit
\eqn{A(R)}{R^i{}_a{}^j{}_b t^a{}_k t^b{}_l=t^j{}_b t^i{}_a R^a{}_k{}^b{}_l,\
\Delta t^i{}_j=t^i{}_a\tens t^a{}_j,\ \eps t^i{}_j=\delta^i{}_j,\ {\rm i.e.}\
R\vect_1\vect_2=\vect_2\vect_1R, \ \Delta \vect=\vect\tens\vect,\
\eps\vect=\id}
where $\vect_1$ and $\vect_2$ denote $\vect\tens\id$ and $\id\tens \vect$ in
$M_n\tens M_n$ with values in $A(R)$.

The power of this matrix notation lies in the fact that $\tens$ is used only
to
refer to the abstract tensor
product of copies of the algebra (as in defining the axioms of a Hopf algebra
etc). The matrix tensor product as in $M_n\tens M_n$ is suppressed and its
role
is replaced by the suffices ${}_1,{}_2$ etc when needed. Thus $R=R_{12}$ (with
the indices suppressed when there are only two $M_n$ in the picture) while
$\vect_1$ means that the $\vect$ is viewed as a matrix in the first $M_n$
(with
values in $A(R)$). The rules of the notation are that matrices are
understood as multiplied in the usual order, independently in the ${}_1,{}_2$
etc copies of $M_n$. Using this notation (or directly with indices) one can
see
at once that $A(R)$ has two fundamental representations $\rho^\pm$ in $M_n$
defined by
\eqn{funrep}{ \rho^+(t^i{}_j)^k{}_l=R^i{}_j{}^k{}_l, \quad
\rho^-(t^i{}_j)^k{}_l=R^{-1}{}^k{}_l{}^i{}_j,\quad {\rm i.e.}\quad
\rho^+_2(\vect_1)=R_{12},\quad \rho^-_2(\vect_1)=R_{21}^{-1}.}
In the compact notation the proof reads
\align{\rho^+_3(R_{12}\vect_1\vect_2)&&=R_{12}\rho_3^+(\vect_1)
\rho_3^+(\vect_2)=R_{12}R_{13}R_{23}\\
&&=R_{23}R_{13}R_{12}=\rho^+_3(\vect_2)\rho_3^+(\vect_1)R_{12}
=\rho^+_3(\vect_2\vect_1 R_{12})\\
\rho^-_3(R_{12}\vect_1\vect_2)&&=R_{12}\rho_3^-(\vect_1)\rho_3^-(\vect_2)
=R_{12}R_{31}^{-1}R_{32}^{-1}\\
&&=R_{32}^{-1}R_{31}^{-1}R_{12}=\rho^-_3(\vect_2)\rho_3^-(\vect_1)R_{12}
=\rho^-_3(\vect_2\vect_1 R_{12}).}
Note that for $\rho^-$ we need
$R_{12}R^{-1}_{31}R^{-1}_{32}=R^{-1}_{32}R^{-1}_{31}R_{12}$, i.e.
$R_{31}R_{32}R_{12}=R_{12}R_{32}R_{31}$ which is again the QYBE after a
relabeling of the positions in $M_n^{\tens3}$.

This bialgebra $A(R)$ is important because for the standard $R$-matrices one
has a convenient construction of the quantum function algebras $\CO_q(G)$
deforming the ring of representative functions on compact simple group
$G$\cite{FRT:lie}. One has to quotient the bialgebra by suitable further
relations (or localise a determinant) to obtain a Hopf algebra. In the
standard
case of course one knew that the result was dual quasitriangular because of
Drinfeld's result that $U_q(g)$ was quasitriangular (over formal
power-series).
The question for a general $R$-matrix was not so clear at the time and was
resolved in \cite{Ma:qua}\cite{Ma:seq} where we showed that there is always
some kind of quasitriangular structure in the form of a map $\CR:A(R)\to
{A(R)^*}^*$. In the more modern setting our result reads as follows. One can
also see subsequent works such as \cite{LarTow:two} but I retain here the
strategy (which comes from physics) of my original proofs, this time with full
pedagogical details.

\begin{propos}  Let $R$ be an invertible solution of the QYBE in $M_n\tens
M_n$. Then the associated matrix bialgebra $A(R)$ is dual quasitriangular with
$\CR:A(R)\tens A(R)\to k$ given by
$\CR(\vect\tens 1)=\id=\CR(1\tens \vect)$ and $\CR(\vect_1\tens\vect_2)=R$
extended as a bialgebra bicharacter according to (\ref{dqua1}). Explicitly,
\align{\CR(t^{i_1}{}_{j_1}t^{i_2}{}_{j_2}\cdots
t^{i_M}{}_{j_M}\tens t^{k_N}{}_{l_N} t^{k_{N-1}}{}_{l_{N-1}}\cdots
t^{k_1}{}_{l_1})\nqquad\nqquad\nqquad\nqquad\nqquad\nqquad\nqquad
\nqquad\nquad&&\\
&&\nquad=R^{i_1}{}_{m_{11}}{}^{k_1}{}_{n_{11}}
\ \, R^{m_{11}}{}_{m_{12}}{}^{k_2}{}_{n_{21}}\ \
\cdots \ R^{m_{1N-1}}{}_{j_{1}}{}^{k_N}{}_{n_{N1}}\\
&&\ \, R^{i_2}{}_{m_{21}}{}^{n_{11}}{}_{n_{12}}
R^{m_{21}}{}_{m_{22}}{}^{n_{21}}{}_{n_{22}}\\
&&\qquad\ \vdots\qquad\qquad\qquad\qquad\qquad\qquad\quad\vdots\cr
&&\ R^{i_M}{}_{m_{M1}}{}^{n_{1M-1}}{}_{l_1} \ \ \cdots\quad \ \, \cdots\quad
R^{m_{MN-1}}{}_{j_{M}}{}^{n_{NM-1}}{}_{l_N}=Z_R({\scriptstyle K\atop
{{{\scriptstyle I
\lform
J}\atop {\scriptstyle L}}}})}
where the last notation is as a partition function\cite[Sec. 5.2]{Ma:qua}.
Here
$I=(i_1,\cdots ,i_M)$ and $K=(k_1,\cdots k_N)$ etc. There is a similar
expression for $\CR^{-1}$. If we  adopt the notation $\bar K=(k_N,\cdots k_1)$
and $t^{i_1}{}_{j_1}\cdots t^{i_M}{}_{j_M}=t^I{}_J$ then
\[ \CR(t^I{}_J\tens t^{\bar K}{}_{\bar L})=Z_R({\scriptstyle K\atop
{{{\scriptstyle I
\lform
J}\atop {\scriptstyle L}}}}),\quad \CR^{-1}(t^{\bar I}{}_{\bar J}\tens
t^K{}_L)=Z_{R^{-1}}({\scriptstyle K\atop {{{\scriptstyle I
\lform
J}\atop {\scriptstyle L}}}})\]
\end{propos}
\proof Note that $\CR(\vect_1\tens\vect_2)=\rho^+_2(\vect_1)$ and the proof
that this extends in its first input  as a bialgebra bicharacter is exactly
the
proof above that $\rho^+$ extends to products as a representation. Thus we
have
$\CR(a\tens\vect)=\rho^+(a)$ for all $a$ and
$\CR(ab\tens\vect)=\CR(a\tens\vect)\CR(b\tens\vect)$ as we require for the
first of (\ref{dqua1}). In particular,
\[\CR(\vect_1\vect_2\cdots\vect_M\tens\vect_{M+1})=\rho^+_{M+1}
(\vect_1\vect_2\cdots\vect_M)=R_{1 M+1}\cdots R_{MM+1}\]
is well-defined. Next the tensor product of representations is also a
representation (because $A(R)$ is a bialgebra), hence there is a well-defined
algebra map $\rho^+{}^{\tens N}:A(R)\to M_n^{\tens N}$ given by
$\rho^+{}^{\tens
N}(a)=(\rho^+_1\tens\rho^+_2\tens\cdots\tens\rho^+_N)\circ\Delta^{N-1}(a)$. In
particular,
\align{{}\nqquad&&R_{1M+1 }R_{1M+2}\cdots R_{1M+N}\\
&&R_{2M+1}R_{2M+2}\cdots R_{2M+N}=\rho^+_{M+1}(\vect_1\cdots\vect_M)\cdots
\rho^+_{M+N}(\vect_1\cdots\vect_M)=(\rho^+)^{\tens N}(\vect_1\cdots\vect_M)\\
&&\quad \cdots\qquad\cdots \\
&&R_{MM+1}\cdots\quad \cdots R_{MM+N}}
also depends only on $\vect_1\vect_2\cdots\vect_M$ as an element of $A(R)$.
The
array on the left can be read (and multiplied up) column after column (so that
the first equality is clear) or row after row (like reading a book). The two
are the same when we bear in mind that $R$ living in distinct copies of
$M_n\tens M_n$ commute. The expression  is just the array $Z_R$ in our compact
notation. If we define
$\CR(\vect_1\vect_2\cdots\vect_M\tens\vect_{M+N}\cdots\vect_{M+1})$ as this
array, we know that the second of (\ref{dqua1}) will hold and that $\CR$ is
well defined in its first input.

Now we repeat the steps above for the second input of $\CR$. Thus
$\CR(\vect_1\tens\vect_2)=R=\bar\rho^+_1(\vect_2)$ extends in its second input
as a bialgebra bicharacter since this ${\bar\rho}^+$ extends as an
antirepresentation $A(R)\to M_n$ (proof similar to that for $\rho^+$). Thus we
define $\CR(\vect\tens a)={\bar\rho}^+(a)$ and
in particular,
\[ \CR(\vect_M\tens
\vect_{M+N}\cdots\vect_{M+2}\vect_{M+1})={\bar\rho}^+_M(\vect_{M+N}
\cdots\vect_{M+2}\vect_{M+1})=R_{MM+1} R_{M M+2} \cdots R_{MM+N}\]
is well-defined. Likewise, we can take tensor products of ${\bar\rho}^+$ and
will again have well-defined anti-representations. Hence
\align{&&\nqquad R_{1M+1}R_{1M+2}\cdots R_{1M+N}\\
&&\nqquad R_{2M+1}R_{2M+2}\cdots
R_{2M+N}={\bar\rho}^+_1(\vect_{M+N}\cdots\vect_{M+1})\cdots
{\bar\rho}^+_{M}(\vect_{M+N}\cdots\vect_{M+1})=({\bar\rho}^+)^{\tens
M}(\vect_{M+N}\cdots\vect_{M+1})\\
&&\nqquad \quad \cdots\qquad\cdots\\
&&\nqquad R_{MM+1}\cdots\quad \cdots R_{MM+N}}
depends only on $\vect_{M+N}\cdots\vect_{M+2}\vect_{M+1}$ as an element of
$A(R)$. The first equality comes from writing out the ${\bar\rho}^+$ and
rearranging the resulting array (bearing in mind that copies of $R$ in
distinct
$M_n$ tensor factors commute). The resulting array of matrices then coincides
with that above, which we have already defined as
$\CR(\vect_1\cdots\vect_M\tens\vect_{M+N}\cdots\vect_{M+1})$. We see that this
array then is well defined as a map $A(R)\tens A(R)\to k$, in its fist input
(for fixed $\vect_{M+N},\cdots,\vect_{M+1}$) by its realization as a tensor
power of $\rho^+$ and in its second input (for fixed $\vect_1,\cdots,\vect_M$)
by its realization as a tensor power of ${\bar\rho}^+$. By its construction,
it
obeys (\ref{dqua1}).

Next, we note that when $R$ is invertible there is a similar construction for
$\CR^{-1}$ to that for $\CR$ above. Here $\CR^{-1}$ obeys equations similar to
(\ref{dqua1}) but with its second input multiplicative and its first input
antimultiplicative. We use $R^{-1}$ in the role of $R$, for example,
$\CR^{-1}(\vect\tens a)=\rho^-(a)$ extends as a
representation, while $\CR^{-1}(a\tens\vect)$ extends as an
antirepresentation.
The steps are entirely analogous to those above, and we arrive at the
partition
function $Z_{R^{-1}}$. We have to show that $\CR,\CR^{-1}$ are inverse in the
convolution algebra of maps $A(R)\tens A(R)\to k$. Explicitly, we need,
\eqn{multiRinv}{\CR(t^I{}_A\tens t^{\bar K}{}_{\bar B})\CR^{-1}(t^{\bar
A}{}_{\bar J}\tens t^B{}_L)=\delta^I{}_J\delta^K{}_L,\quad{\rm i.e.}\quad
Z_R({\scriptstyle \bar K\atop {{{\scriptstyle I
\lform
A}\atop {\scriptstyle \bar B}}}})Z_{R^{-1}}({\scriptstyle  B\atop
{{{\scriptstyle \bar A
\lform
\bar J}\atop {\scriptstyle L}}}})=\delta^I{}_J\delta^K{}_L}
and similarly on the other side. Writing the arrays in our compact notation we
have
\align{&&R_{1 M+N}\cdots R_{1M+1}\\
&&\cdots\qquad\cdots\\
&&R_{M-1 M+N}\cdots R_{M-1 M+1}\\
&&R_{M M+N}\cdots R_{MM+1}R^{-1}_{MM+1}\cdots R_{MM+N}^{-1}\\
&&\phantom{R_{M M+N}\cdots R_{MM+N}}R^{-1}_{M-1M+1}\cdots R^{-1}_{M-1M+N}\\
&&\phantom{R_{M M+N}\cdots R_{MM+N}}\cdots\qquad\cdots\\
&&\phantom{R_{M M+N}\cdots R_{MM+N}}R^{-1}_{1M+1}\cdots R^{-1}_{1M+N}}
Here the copies of $M_n$ numbered $1\cdots M$ on the left correspond to the
index $I$, on the right to $\bar J$ (they occur reversed). The copies of $M_n$
numbered $M+1\cdots M+N$ correspond on the top to $\bar K$ (they occur
reversed) and on the bottom to $L$. In between they are matrix-multiplied as
indicated, corresponding to the sum over $A,B$. The overlapping line here
collapses after cancellation of inverses ending in $\id$ in the copy of $M_n$
numbered $M$, and results in a similar picture with one row less. Repeating
this, the whole thing collapses
to the identity in all the copies of $M_n$.

Finally, we check (\ref{dqua2}), which now takes the form
\eqn{multiA(R)}{ t^{\bar K}{}_{\bar B} t^I{}_A Z_R({\scriptstyle B\atop
{{{\scriptstyle A\lform J}\atop {\scriptstyle L}}}})
= Z_R({\scriptstyle K\atop {{{\scriptstyle I\lform A}\atop {\scriptstyle
B}}}})t^A{}_J t^{\bar B}{}_{\bar L}.}
In the compact notation we compute
\align{&&\vect_{M+N}\cdots\vect_{M+1}\vect_1\cdots\vect_M R_{1M+1}\cdots
R_{1M+N}\\
&&\phantom{\vect_{M+N}\cdots\vect_{M+1}\vect_1\cdots\vect_M}\cdots
\qquad\cdots\\
&&\phantom{\vect_{M+N}\cdots\vect_{M+1}\vect_1\cdots\vect_M}R_{MM+1}\cdots
R_{MM+N}\\
&&=R_{1M+1}\cdots
R_{1M+N}\vect_1\vect_{M+N}\cdots\vect_{M+1}\vect_2\cdots\vect_M R_{2M+1}\cdots
R_{2M+N}\\
&&\phantom{=R_{1M+1}\cdots
R_{1M+N}\vect_1\vect_{M+N}\cdots\vect_{M+1}\vect_2\cdots\vect_M}
\cdots\qquad\cdots\\
&&\phantom{=R_{1M+1}\cdots
R_{1M+N}\vect_1\vect_{M+N}\cdots\vect_{M+1}\vect_2\cdots\vect_M}R_{MM+1}\cdots
R_{MM+N}\\
&&=\\
&&\vdots\\
&&=R_{1M+1 }R_{1M+2}\cdots R_{1M+N}\\
&&\phantom{=}R_{2M+1}R_{2M+2}\cdots R_{2M+N}\\
&&\phantom{=}\quad \cdots\qquad\cdots \\
&&\phantom{=}R_{MM+1}\cdots\quad \cdots R_{MM+N}
\vect_1\vect_2\cdots\vect_M\vect_{M+N}\cdots\vect_{M+1}.}
Here the copies of $M_n$ numbered $1,\cdots, M$ on the left correspond to the
index $I$, and the copies of $M_n$ numbered $M+1,\cdots ,M+N$ correspond on
the
top to $\bar K$ (they occur reversed), etc. The first equality makes repeated
use of the relations (\ref{A(R)}) of $A(R)$ to give
\[   \vect_{M+N}\cdots\vect_{M+1}\vect_1R_{1M+1}\cdots R_{1M+N}=R_{1M+1}\cdots
R_{1M+N}\, \vect_1 \vect_{M+N}\cdots\vect_{M+1}.\]
The $\vect_2\cdots\vect_M$ move past the  $R_{1M+1}\cdots R_{1M+N}$ freely
since they live in different matrix spaces. This argument for the first
equality is then applied to move $\vect_{M+N}\cdots\vect_{M+1}\vect_2$, and so
on.
The arguments in this proof may appear complicated, but in fact this kind of
repeated matrix multiplication (multiplication of entire rows or columns of
matrices) is quite routine in the context of exactly solvable statistical
mechanics (where the QYBE originated).  \endproof

Let us note that while the algebra relations (\ref{A(R)}) of $A(R)$ do not
depend on the normalization of $R$, the dual quasitriangular structure does.
The elements $t^{i_1}{}_{j_1}\cdots t^{i_M}{}_{j_M}$ of $A(R)$ have a
well-defined degree $|t^{i_1}{}_{j_1}\cdots t^{i_M}{}_{j_M}|=M$ (the algebra
is
graded), and if $R'=\lambda R$ is a non-zero rescaling of our solution $R$
then
the corresponding dual quasitriangular structure is changed to
\eqn{scalingR}{\CR'(a\tens b)=\lambda^{|a||b|}\CR(a\tens b)}
 on homogeneous elements. This is evident from the expression in terms of
$Z_R$
that we have obtained in the last proposition. Finally, in view of the reasons
that we passed to the dual setting it is obvious that

\begin{theorem} Let $(A,\CR)$ be a dual quasitriangular bialgebra. Then
$\CM^A$
the category of right $A$-comodules is braided. In the Hopf algebra case the
finite-dimensional comodules are rigid,
\[ \Psi_{V,W}(v\tens w)=\sum w\bo\tens v\bo\CR(v\bt\tens w\bt), \quad
\beta_{V^*}(f)=(f\tens S)\circ\beta_V\]
where $\beta_{V^*}(f)\in V^*\tens A$ is given as a map $V\to A$.
\end{theorem}
\proof This is an entirely trivial dualization of the proof of Theorem~1.10
above. For example, $\Psi$ is an intertwiner because
\align{\Psi_{V,W}(v\tens w)&\mapsto&\sum w\bo\bo\tens v\bo\bo\tens w\bo\bt
v\bo\bt\CR(v\bt\tens w\bt)\\
&&=\sum w\bo\tens v\bo\tens w\bt\o v\bt\o\CR(v\bt\t\tens w\bt\t)\\
&&=\sum w\bo\tens v\bo\tens \CR(v\bt\o\tens w\bt\o) v\bt\t w\bt\t\\
&&=\sum w\bo\bo\tens v\bo\bo\tens \CR(v\bo\bt\tens w\bo\bt) v\bt w\bt}
where the arrow is the tensor product $W\tens V$ coaction. We used
(\ref{dqua2}). The result is $\Psi$ applied to the result of the tensor
product
coaction $V\tens W$. The hexagons (\ref{Psi-hex}) correspond in a similarly
trivial way to (\ref{dqua1}).
\endproof

Note that $\CC(V,R)$ in Example~1.1 forms a subcategory of $\CM^{A(R)}$.
Moreover, in the dualizable case there is a Hopf algebra $GL(R)\supset A(R)$
such that $\CC(V,V^*,R)$ is a subcategory of $\CM^{GL(R)}$. The relevant
coactions are
\eqn{funcorep}{e_i\mapsto e_j\tens t^j{}_i,\quad f^i\mapsto f^j\tens S
t^i{}_j}
and we recover from Theorem~1.16 the braidings quoted. Also, in
\cite{Ma:pro}\cite{Ma:tan}\cite{Ma:eul}\cite{Ma:bg} we regarded this
proposition as a starting point and
set out to prove something further, namely its converse. If $\CM^A$ is braided
then $A$ has induced on it by Tannaka-Krein reconstruction a
dual quasitriangular structure. We will see this in Section~3. \cite{Ma:tan}
generalised the Tannaka-Krein theorem to the setting of dual quasi-Hopf
algebras (associative up to an isomorphism cf\cite{Dri:quas}) while
\cite{Ma:eul}\cite{Ma:bg} generalised it to the braided setting. It is more or
less the {\em sine qua non} for the work here.

\section{Braided Tensor Product Algebra and Braided Hopf Algebras}

So far we have described braided monoidal or quasitensor categories and ways
to
obtain them. Now we
begin our main task and study algebraic structures living in such categories.
For this we use the
diagrammatic notation of Section~1.1. Detailed knowledge of quantum groups etc
is not required in this section.

The idea of an algebra $B$ in a braided category is just the usual one. Thus
there should be product and unit morphisms
\eqn{bra-alg}{\cdot: B\tens B\to B,\qquad\eta:\und 1\to B}
obeying the usual associativity and unity axioms but now as morphisms in the
category. Note that the term `algebra' is being used loosely since we have not
discussed direct sums and linearity under a field or ring. These notions are
perfectly compatible with what follows but do not play any particular role in
our general constructions.

The fundamental lemma for the theory we need is the generalization to this
setting of the usual $\Z_2$-graded or super-tensor product of superalgebras:
\begin{lemma}\cite{Ma:bra}--\cite{Ma:exa} Let $B,C$ be two algebras in a
braided category. There is a {\em
braided tensor product
algebra} $B\und\tens C$, also living in the braided category. It has product
$(\cdot_B\tens\cdot_C)\circ(\id\tens\Psi_{C,B}\tens\id)$ and tensor product
unit morphism.
\end{lemma}
\proof We repeat here the diagrammatic proof\cite{Ma:sta}. The box is the
braided tensor product multiplication,
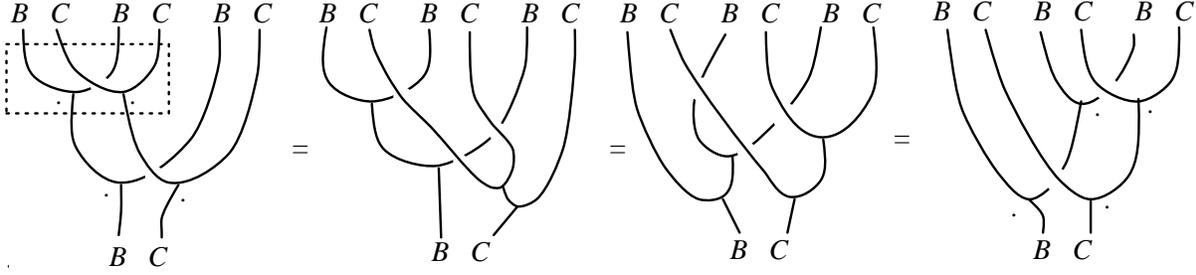
\[ \epsfbox{tensprodproof.eps}\]
The first step uses functoriality as in (\ref{Psi-funct-bra}) to pull the
product morphism through the braid crossing. The second equality uses
associativity of the products in $B,C$ and the third equality uses
functoriality again in reverse. The product is manifestly a morphism in the
category because it is built out of morphisms. Finally, the unit is the tensor
product
one because the braiding is trivial on
$\und 1$. \endproof

In the concrete case the braided tensor product is generated by $B=B\tens 1$
and $C=1\tens C$ and an exchange law between the two factors given by $\Psi$.
This is because $(b\tens 1)(1\tens c)$ $=$ $(b\tens c)$ while $(1\tens
c)(b\tens
1)=\Psi(c\tens b)$. Another notation is to label the elements of the second
copy in the braided tensor product by ${}'$. Thus $b\equiv (b\tens 1)$ and
$c'\equiv (1\tens c)$. Then if $\Psi(c\tens b)=\sum b_k \tens c_k$ say, we
have
the braided-tensor product  relations
\eqn{btens-stat}{ c'b\equiv (1\tens c)(b\tens 1)=\Psi(c\tens b)=\sum b_k\tens
c_k\equiv \sum b_k c'_k.}
This makes clear why the lemma generalizes the notion of
$\Z_2$-graded or super-tensor product. Note also that there is an equally good
{\em opposite braided tensor product}
with the inverse braid crossing in Lemma~2.1. This is simply the braided
tensor
product
algebra constructed in the mirror-reversed category $\bar{\CC}$ but with the
result viewed in our original category.

\subsection{Braided Hopf Algebras}

Armed with the braided tensor product of algebras in a braided category we can
formulate the notion of
Hopf algebra.

\begin{defin}\cite{Ma:bra}--\cite{Ma:exa} A Hopf algebra in a braided
category or {\em braided-Hopf algebra} is $(B,\Delta,\eps,S)$ where $B$ is an
algebra in the category and $\Delta:B\to B\und\tens B$,
$\eps:B\to\und 1$ are algebra homomorphisms where $B\und\tens B$ has the
braided tensor product algebra structure. In addition, $\Delta,\eps$ obey the
usual
coassociativity and counity axioms to form a coalgebra in the
category, and $S:B\to B$ obeys the usual axioms of an antipode.
If there is no antipode then we speak of a braided-bialgebra or bialgebra in a
braided category.
\end{defin}

In diagrammatic form the algebra homomorphism and braided-antipode axioms read
\eqn{hopf-ax}{\epsfbox{hopf-ax.eps}}
One then proceeds to develop the usual elementary theory for these Hopf
algebras.
For example, recall that the usual antipode is an antialgebra map.

\begin{lemma} For a braided-Hopf algebra $B$, the braided-antipode obeys
$S(b\cdot c)=\cdot\Psi(Sb\tens Sc)$ and $S(1)=1$, or more abstractly,
$S\circ\cdot=\cdot\circ\Psi_{B,B}\circ(S\tens S)$ and $S\circ\eta=\eta$.
\end{lemma}
\proof In diagrammatic form the proof is\cite{Ma:tra}
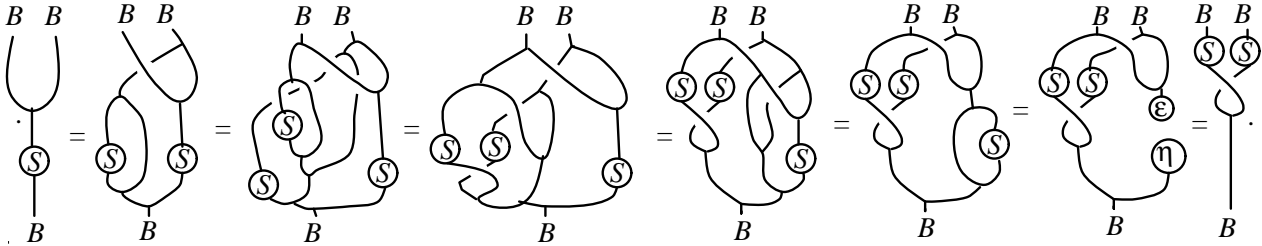
\[\epsfbox{antproof.eps}\]
In the first two equalities we have grafted on some circles containing the
antipode, knowing they are trivial from (\ref{hopf-ax}). We then use the
coherence theorem to lift the second $S$ over to the left, and associativity
and coassociativity to reorganise the branches. The fifth equality uses the
axioms (\ref{hopf-ax}) for $\Delta$.
\endproof

Here we want to mention a powerful {\em input-output symmetry} of these
axioms.
Namely, turn the
pages of
this book upside down and look again at these diagrams. The axioms of a
braided-Hopf algebra (\ref{hopf-ax})
are unchanged except that the roles of product/coproduct and unit/counit
morphisms are interchanged. The proof of Lemma~2.1 becomes the proof of a new
lemma expressing coassociativity of the braided-tensor product of two
coalgebras. Meanwhile the proof of Lemma~2.3 reads as the proof of a new lemma
that the braided-antipode is a braided-anti-coalgebra map.

This applies therefore to all results that we prove about bialgebras or Hopf
algebras in braided categories provided all notions are suitably turned
up-side-down. This is completely rigorous and nothing to do with
finite-dimensionality or individual dual objects. In addition, there is a {\em
left-right symmetry} of the axioms
consisting of reflecting in mirror about a vertical axis combined with
reversal
of all braid crossings. These symmetries of the axioms can be taken together
so
that we obtain precisely four theorems for the price of one when we use the
diagrammatic method.

An endemic problem for those working in Hopf algebras is that every time
something is proven one has to laboriously figure out its
input-output-reversed version or its version with opposite left-right
conventions. This problem is entirely solved by reflecting in a mirror or
turning up-side-down.

\subsection{Dual Braided Hopf Algebras}

Suppose now that the category has dual objects (is rigid) in the sense
explained in Section~1.1. In this case the input-output symmetry of the axioms
of a Hopf algebra becomes realised concretely as the construction of a dual
Hopf algebra.

\begin{propos} If $B$ is a braided-Hopf algebra, then its left-dual $B^*$ is
also a braided-Hopf algebra with product, coproduct, antipode, counit and unit
given by
\[\epsfbox{B-dual.eps}\]
\end{propos}
\proof Associativity and coassociativity follow at once from coassociativity
and associativity of $B$. Their crucial compatibility property comes out as
\[ \epsfbox{dual-bialg.eps}\]
where we use the double-bend axiom (\ref{lduals-bra}) for dual objects. The
antipode property comes out just as easily.
\endproof

We see that in diagrammatic form the dual bialgebra or Hopf algebra is
obtained
by rotating the desired structure map in an anticlockwise motion and without
cutting any of the attaching strings. For right duals the motion should be
clockwise. Let us stress that this dual-Hopf algebra construction of an
individual object should not be confused with the rather more powerful
input-output symmetry for the axioms introduced above.

\subsection{Braided Actions and Coactions}

Another routine construction is the notion of module and its
input-output-reversed notion of comodule. These are just the obvious ones but
now as morphisms in the category. Let us check that the tensor product of
modules of a bialgebra is a module. If $V,\alpha_V:B\tens V\to V$ and
$W,\alpha_W:B\tens W\to W$ are two left modules then
\eqn{tensVW}{\epsfbox{tensVWproof.eps}}
is the definition (in the box) of tensor product module and proof that it is
indeed a module. The first equality is the homomorphism property of $\Delta$.
Likewise for right modules by left-right reflecting the proof in a mirror and
also reversing all braid crossings.

As for usual Hopf algebras, the left-right symmetry can be concretely realised
via the antipode. Thus if $V,\alpha^R:V\tens B\to V$ is a right module then
\eqn{V-rlmod}{\epsfbox{V-rlmodproof.eps}}
shows the construction of the corresponding left module. This is shown in the
box.

Also if the left-module $V$ has a left-dual $V^*$ then this becomes a
right-module with
\eqn{V*-mod}{\epsfbox{Vs-mod.eps}.}
Combining these two constructions we conclude (with the obvious definition of
intertwiners or morphisms between braided modules):

\begin{propos}\cite{Ma:bg} Let $B$ be a bialgebra in the  braided category
$\CC$. Then the category ${}_B\CC$ of braided left-modules is a monoidal
category. If $B$ is a Hopf algebra in $\CC$ and $\CC$ is rigid then ${}_B\CC$
is rigid.
\end{propos}
\proof For the second part we feed the result of (\ref{V*-mod}) into
(\ref{V-rlmod}). A more traditional-style proof with commuting
diagrams is in \cite{Ma:bg} for comparison. It should convince the reader of
the power of the diagrammatic method. \endproof

Naturally, a braided left $B$-module algebra is by definition an algebra
living
in the category ${}_B\CC$. This means an algebra $C$ such that
\eqn{C-modalg}{\epsfbox{C-modalg.eps}\qquad B{\rm -Module\  Algebra}}
Likewise for other constructions familiar for actions of bialgebras or Hopf
algebras. For example, a coalgebra $C$ in the category ${}_B\CC$ is a
coalgebra
such that
\eqn{C-modcoalg}{\epsfbox{C-modcoalg.eps}\qquad B{\rm -Module\ Coalgebra}.}

For the right-handed theory reflect the above in a mirror and reverse all
braid
crossings. This gives the notion of right $B$-module algebras etc. Next, by
turning the pages of this book upside down we have all the
corresponding results for comodules in place of modules. Thus the category
$\CC^B$ of right-comodules had a tensor product and in the Hopf algebra and
rigid case is also rigid. Likewise for left comodules.

\begin{propos}\cite{Ma:bos} If $C$ is a left $B$-module algebra then there is
a
braided {\em cross product} or semidirect product algebra $C\cocross B$ built
on the object $C\tens B$. Likewise for right $B$-module algebras, left
$B$-module coalgebras and right $B$-comodule coalgebras. The semidirect
(co)product maps are
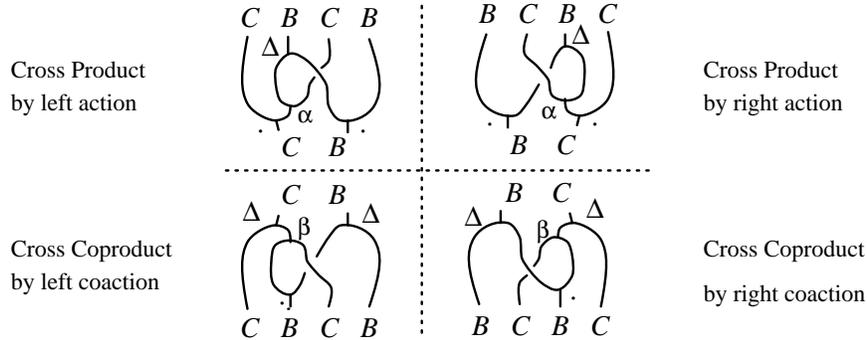
\[ \epsfbox{C-cross.eps}\]
\end{propos}
\proof We only need to prove one of these by our diagrammatic means to
conclude
all four. Full details are in \cite{Ma:bos}. \endproof

\begin{example}\cite[Appendix]{Ma:lin} Let $B$ be a braided-Hopf algebra. Then
$B$ is a left $B$-module algebra by the {\em braided adjoint} action
$\cdot^2\circ(\id\tens\Psi_{B,B})\circ(\id\tens S\tens\id)\circ
(\Delta\tens\id)$.
\end{example}

\begin{example}\cite{Ma:lie} Let $B$ be a braided-Hopf algebra with left dual
$B^*$. Then $B$ is a right $B^*$-module algebra by a {\em braided right
regular} action $\ev_B\circ(\id\tens
S\tens\id)\circ(\id\tens\Delta)\circ\Psi_{B,B^*}$
\end{example}

The verification of these examples is a nice demonstration of the techniques
above. In diagrammatic form they read
\eqn{Ad-rvect}{\epsfbox{Adfrag.eps}\qquad\quad\epsfbox{regfrag.eps}}
The adjoint action leads to a notion of braided Lie algebra\cite{Ma:lie} among
other applications, while the right regular action corresponds to the action
of
fundamental vector fields. It can also be used to construct a braided Weyl
algebra cf \cite{Ma:fre}.

\subsection{Braided-(Co)-Commutativity}

Next we come to the question of commutativity or cocommutativity in a braided
category. Again, we only have to work
with one of these and turn our diagrams up-side-down for the other. The main
problem is that the naive opposite-coproduct
\eqn{opcoprod}{\bar\Delta=\Psi^{-1}_{B,B}\circ\Delta}
does {\em not} make $B$ into a bialgebra in our original braided category
$\CC$, but rather gives a bialgebra in the mirror-reversed category
$\bar{\CC}$. Thus there is a notion of opposite bialgebra (and if $B$ is a
Hopf
algebra with invertible antipode then $S^{-1}$ provides an antipode) but it
forces us to leave the category.

Hence there is no way to consider bialgebras or Hopf algebras that are
cocommutative in the sense that they
coincide with their opposite. There is so far no intrinsic notion of
braided-cocommutative Hopf algebras for this
reason. On the other hand we have introduced in \cite{Ma:bra} the notion of a
braided-cocommutativity with respect to a module. This is a property of a
module on which $B$ acts.

\begin{defin}\cite{Ma:bra} A braided left module $(V,\alpha_V)$ is {\em
braided-cocommutative} (or $B$ is braided-cocommutative with respect to $V$)
if
\[\epsfbox{V-cocom.eps}\qquad {\rm Braided-Cocommutativity}\]
\end{defin}

To understand this notion suppose that the category is symmetric not strictly
braided. In this case $\Psi^2=\id$ and we see that the condition is implied by
$\bar\Delta=\Delta$. But in a general braided category we cannot disentangle
$V$ and must work with this weaker notion. Moreover, as far as such modules
are concerned the bialgebra $B$ has all the usual representation-theoretic
features of usual cocommutative Hopf algebras. One of these is that their
tensor product is symmetric under the usual transposition of the underlying
vector spaces of modules. The parallel of this is

\begin{propos}\cite{Ma:tra} Let $B$ be a bialgebra in a braided category and
define $\CO(B)\subset {}_B\CC$ the subcategory of braided-cocommutative
modules. Then $\CO(B)$ is closed under $\tens$. Moreover, the tensor product
in
$\CO$ is braided with braiding induced by the braiding in $\CC$,
\[\epsfbox{comtensbox.eps}\qquad {\rm Braided-Commutativity\ of\ Product\
of\ Modules}\]
\end{propos}
\proof The first part is given in detail in \cite{Ma:tra} in a slightly more
general context.
The second part follows from Definition~2.9 by adding an action on $W$ to both
sides. \endproof

The trivial representation is always braided-cocommutative. In many examples
the adjoint representation in Example~2.7 is also braided-cocommutative. In
this case one can formulate properties like those of an enveloping algebra of
a braided-Lie algebra\cite{Ma:lie}. One can formally define a {\em braided
group} as a pair consisting of a Hopf algebra in a braided category and a
class
of braided-cocommutative modules. This turns out to be a useful notion because
in many situations it is only this weak notion of cocommutativity that is
needed. For example

\begin{theorem}\cite{Ma:bos} Let $B,C$ be Hopf algebras in a braided category
and $C$ a braided-cocommutative-$B$-module algebra and coalgebra. Then
$C\cocross B$ forms
a Hopf algebra in the braided category with the braided tensor product
coalgebra structure.
\end{theorem}

\subsection{Quantum-Braided Groups}

We can go further and consider bialgebras that are quasi-cocommutative in
some sense, analogous to the idea of a quasitriangular bialgebra in
Section~1.2. To do this we require a second coproduct which we denote
$\Deltaop: B\to B\tens B$ also making $B$ into a bialgebra. In this case we
have the notion of a braided $B$-module with respect to which $\Deltaop$
behaves like an opposite coproduct. This is just as in Definition~2.9 but with
the left hand $\Delta$ replaced by $\Deltaop$. The class of such $B$-modules
is
denoted $\CO(B,\Deltaop)$.

The second ingredient that we need is a quasitriangular structure which is
understood now as a convolution-invertible morphism $\CR:\und 1\to B\tens B$.
With these ingredients the analogue of Definition~1.5 is\cite{Ma:tra}
\eqn{B-univR}{\epsfbox{B-univR.eps}.}
The braided analogue of Theorem~1.10 is then

\begin{theorem}\cite{Ma:tra} Let $(B,\Deltaop,\CR)$ be a quasitriangular
bialgebra in a braided category. Then $\CO(B,\Deltaop)\subset {}_B\CC$ is a
braided monoidal category with braiding
$\Psi^{\CO}_{V,W}=\Psi_{V,W}\circ(\alpha_V\tens\alpha_W)
\circ\Psi_{B,V}\circ(\CR\tens\id)$.
\end{theorem}
\proof Diagrammatic proofs are in \cite[Sec. 3]{Ma:tra}. \endproof

The dual theory with comodules and an opposite product is developed in
\cite{Ma:eul}\cite{Ma:bg}. We mention here only that turning (\ref{B-univR})
up-side-down and then setting the category to be the usual one of vector
spaces
returns not the axioms of a dual quasitriangular structure as in Section~1.3
but its inverse. This reversal is due to the fact that the categorical
dualization in Section~2.2 yields in the vector space category the opposite
coproduct and product to the usual dualization.

\section{Reconstruction Theorem}

In this section we give a construction (not the only one) for bialgebras and
Hopf algebras in braided categories. There is such a bialgebra associated to a
pair of
braided categories $\CC\to \CV$ or even to a single braided category $\CC$.
The
idea behind this is the theory of Tannaka-Krein reconstruction
generalised to the braided setting.

The Tannaka-Krein reconstruction theorems should be viewed as a generalization
of the simple notion of Fourier Transform. The idea is that the right notion
of
representation of an algebraic structure should itself have enough
structure to reconstruct the original algebraic object. On the other hand many
constructions may appear very simple in terms of the representation theory and
highly non-trivial in terms of the original algebraic object, and vice versa.

In the present setting we know that quantum groups give rise to braided
categories as their representations, while conversely we will see that the
representations or endomorphisms of a category $\CC$ in a category $\CV$ gives
rise to a quantum group in $\CV$.

\subsection{Usual Tannaka-Krein Theorem}

The usual Tannaka-Krein theorem for Hopf algebras says that a monoidal
category
$\CC$ equipped with a functor to the category of vector spaces (i.e. whose
objects can be identified in a strict way with vector spaces) is equivalent to
that of the comodules of a certain bialgebra $A$ reconstructed from $\CC$. All
our
categories $\CC$ are assumed equivalent to small ones.

\begin{theorem} Let $F:\CC\to{\rm Vec}$ be a monoidal functor to the category
of vector spaces with finite-dimensional image. Then there exists a bialgebra
$A$ uniquely
determined as universal with the property that $F$ factors through $\CM^A$. If
$\CC$ is braided then $A$ is dual quasitriangular. If $\CC$ is rigid then $A$
has an antipode. \end{theorem}
\proof We defer this to Theorem~3.11 below. Just set $\CV={\rm Vec}$ there.
\endproof

An early treatment of the bialgebra case is in \cite{Sav:cat}. See also
\cite{DelMil:tan}.
The part concerning the antipode was shown in \cite{Ulb:hop}. That a symmetric
category
gives a dual-triangular structure was pointed out in the modules setting in
\cite{Dri}. See also \cite{Lyu:hop}. It is a trivial step to go from there to
the braided case in which case the result is dual quasitriangular. This has
been done by the author, while at the same time (in order to say something
new)
 generalising in two directions. One is to the quasi-Hopf algebra
setting\cite{Ma:qua}\cite{Ma:tan} and the other to the braided-Hopf algebra
setting\cite{Ma:rec}\cite{Ma:bg}.

This theorem tells us that (dual)quasitriangular Hopf algebras are rather more
prevalent in mathematics (and physics) than we might have otherwise suspected.
It also gives us a useful perspective on any Hopf algebra construction, by
allowing us to go backwards and forwards between representations and the
algebra itself. For example, if we are already given a bialgebra $A$ then
coming out of the reconstruction theorem one has associated to any subcategory
\eqn{CO}{ \CO\subset {}\CM^A}
closed under tensor product, a sub-bialgebra
\eqn{AO}{ A_{\CO}=\cup_{(V,\beta_V)\in \CO} \image(\beta_V)\subset A,\qquad
\image(\beta_V)=\{(f\tens\id)\circ\beta_V(v);\\ v\in V,\ f\in V^*\}.}
If the sub-category is braided then $A_{\CO}$ is dual quasitriangular etc. So
this is a concrete form of the reconstruction theorem in the case where $\CO$
is already in the context of a bialgebra.

For example if $A=A(R)$ and $\CO=\CC(V,R)$ in Section~1 then $A_{\CO}=A(R)$
again. This is because the image of tensor powers of $V$ for the coaction in
(\ref{funcorep}) is clearly any monomial in the generators $\vect$ of $A(R)$.
Hence in this case the subcategory reconstructs all of $A(R)$. The result is
due to Lyubashenko though the proof in \cite{Lyu:hop} is different (and stated
in the triangular case).

For another example let $A$ be a bialgebra and $\CO$ the category of comodules
which are commutative in the sense
$\sum v\bo\tens a v\bt=\sum v\bo\tens v\bt a$ for all $a\in A$ and $v$ in the
comodule. Then $A_{\CO}$ is a bialgebra contained in the center of $A$. This
is
therefore a canonically associated `bialgebra centre' construction.

There are analogous results to these for modules. At the level of Theorem~3.1
the module theory is less powerful only if one functors (as usual) into
familiar finite-dimensional vector spaces.

\goodbreak\goodbreak \goodbreak
\newpage
\subsection{Braided Reconstruction Theorem}

In this section we come to the fully-fledged braided Tannaka-Krein-type
reconstruction theorem. We follow for pedagogical reasons the original module
version\cite{Ma:rec}\cite{Ma:tra}, mainly because the comodule version was
already given in complete
detail in \cite{Ma:bg} and we do not want to repeat it. Also, we give here for
the first time a fully diagrammatic proof.

Throughout this section we fix $F:\CC\to \CV$ a monoidal functor between
monoidal categories. At least $\CV$ should be braided. In this case there is
an
induced functor $V\mapsto \Nat(V\tens F,F)$. We suppose that this functor is
representable. So there is an object $B\in \CV$ such that $\Nat(V\tens
F,F)\isom \Hom_{\CV}(V,B)$ by functorial bijections. Let $\{\alpha_X:B\tens
F(X)\to F(X);\ X\in \CC\}$ be the natural transformation corresponding to the
identity morphism $B\to B$. Then using $\alpha$ and the braiding we get an
induced map
\eqn{rep-mod}{ \Hom_{\CV}(V,B^{\tens n})\to \Nat(V\tens F^n,F^n)}
and we assume that these are likewise bijections. This is the {\em
representability assumption for modules} and we assume it in what follows.

\begin{theorem}\cite{Ma:tra} Let $F:\CC\to \CV$ obey the representability
assumption for modules. Then $B$ is a bialgebra in $\CV$, uniquely determined
as universal with the property that $F$ factors through ${}_B\CV$. If $\CC$ is
braided then $B$ is quasitriangular in the braided category with $\CR$ given
by
the ratio of the braidings in $\CC$ and $\CV$. If $\CC$ is rigid
then $B$ has a braided-antipode.
\end{theorem}

We will give the proof in diagrammatic form. The bijections (\ref{rep-mod})
and
the structure maps in the theorem are
characterized by
\eqn{recon}{\epsfbox{recon.eps}}
Here the assumption that $F$ is monoidal means that there are functorial
isomorphisms $F(X\tens Y)\isom F(X)\tens F(Y)$ and in the rigid case
$F(X^*)\isom F(X)^*$. The latter follow from the uniqueness of duals up
to isomorphism (for example one can define $F(X)^*=F(X^*)$ etc. and any other
dual is
isomorphic). These isomorphisms
are used freely and suppressed in the
notation. The solid node $\alpha_{X\tens Y}$ is $\alpha$ on the composite
object $X\tens Y$ but viewed via the first of these isomorphisms as a morphism
$B\tens F(X)\tens F(Y)\to F(X)\tens F(Y)$. Similarly for $\alpha_{X^*}$. In
this way all diagrams refer to morphisms in $\CV$. The unit $\und 1\to B$
corresponds to the identity natural transformation and the counit to
$\alpha_\und 1$. Their proofs are suppressed.

\goodbreak\goodbreak \goodbreak
\newpage

\begin{lemma} The product on $B$ defined in (\ref{recon}) is associative.
\end{lemma}
\proof We use the definition of $\cdot$ twice in terms of its corresponding
natural transformations and then in reverse
\[ \epsfbox{assoc.eps}\]
Hence the natural transformations corresponding to the two morphisms $B\tens
B\tens B\to B$ coincide and we have an algebra in the category. \endproof

\goodbreak\goodbreak
\begin{lemma} The coproduct $\Delta$ on $B$ defined in (\ref{recon}) is
coassociative.
\end{lemma}
\proof We use the definition of $\Delta$ twice and then in reverse, using in
the middle that $F$ is monoidal
and hence compatible with the (suppressed) associativity in the two categories
\[ \epsfbox{coassoc.eps}\]
The key step is the third equality which follows from functoriality of
$\alpha$ under the associativity morphism $X\tens(Y\tens Z)\to (X\tens Y)\tens
Z$
and that $F$ is monoidal. If $F$ is not monoidal but merely multiplicative one
has here a quasi-associative coproduct as explained in
\cite{Ma:qua}\cite{Ma:tan}.
\endproof

\goodbreak\goodbreak
\begin{propos} The product and coproduct in the last two lemmas fit together
to
form a bialgebra in $\CV$.
\end{propos}
\proof We use the definitions of $\cdot$ and $\Delta$
\[ \epsfbox{bialgebra.eps}\]
\endproof

\goodbreak\goodbreak
\begin{lemma} The second coproduct $\Deltaop$ on $B$ defined in (\ref{recon})
is coassociative.
\end{lemma}
\proof This is similar to the proof of Lemma~3.4
\[ \epsfbox{opcoassoc.eps}\]
\endproof

\goodbreak\goodbreak
\begin{propos} The product and the second coproduct fit together to form a
second bialgebra in $\CV$.
\end{propos}
\proof This is similar to the proof of Proposition~3.5
\[ \epsfbox{opbialgebra.eps}\]
\endproof

\goodbreak\goodbreak
\begin{propos} If $\CC$ is rigid then $S$ defined in (\ref{recon}) is an
antipode for the coproduct $\Delta$.
\end{propos}
\proof The first, second and fourth equalities are the definitions of
$\cdot,S,\Delta$. The fifth uses functoriality of $\alpha$ under the
evaluation
$X^*\tens X\to \und 1$
\[ \epsfbox{antipode.eps}.\]
The result is the natural transformation corresponding to $\eta\circ\eps$.
Similarly for the second line using functoriality under the coevaluation
morphism $\und1\to X\tens X^*$.
\endproof

\goodbreak\goodbreak
\begin{propos} If $\CC$ is braided then $\CR$ defined in (\ref{recon}) makes
$B$ into a quasitriangular bialgebra.
\end{propos}
\proof To prove the first of (\ref{B-univR}) we evaluate the definitions and
use $F$ applied to the hexagon identity in $\CC$ for the fifth equality, and
then in
reverse.
\[ \epsfbox{quasia.eps}.\]
The same strategy works for the second of (\ref{B-univR})
\[ \epsfbox{quasib.eps}\]
Finally, to prove the last of (\ref{B-univR}) we use in the third equality the
functoriality of $\alpha$ under the morphism $\Psi_{X,Y}$:
\goodbreak\goodbreak
\vskip -.2in
\[ \epsfbox{quasic.eps}\]
The construction of $\CR^{-1}$ is based in the inverse natural transformation
to that for $\CR$ and the proof that this then is inverse in the convolution
algebra $\und 1\to B\tens B$ is straightforward using the same techniques.
\endproof

Clearly the definition of the product in (\ref{recon}) is such that $\alpha_X$
become modules. So we have a functor $\CC\to {}_B\CV$. The universal property
of $B$ follows easily from its role as representing object for natural
transformations (\ref{rep-mod}).

\begin{corol} If $\CC$ is braided and $F$ is a tensor functor in the sense
that
the braiding of $\CC$ is mapped on to the braiding of $\CV$ then
$\Deltaop=\Delta$, $\CR$ is trivial and $B$ is a braided group
(braided-cocommutative) with respect to the image of the functor $\CC\to
{}_B\CV$. \end{corol}
\proof  This follows at once from the form of $\CR,\Deltaop$ in (\ref{recon}).
\endproof

For example if $F=\id$ (or the canonical functor into a suitable completion of
$\CC$) then to every rigid braided $\CC$ we have an associated braided Hopf
algebra $B$ and a large class of braided-cocommutative modules $\{\alpha_X\}$.
The ratio of the braidings is trivial and this is why we have from this point
of
view some kind of braided group rather than braided quantum group.

Finally, given a monoidal functor $F:\CC\to \CV$ we can equally well require
representability of the functor
 $V\mapsto \Nat(F,F\tens V)$ and its higher order products i.e., bijections
\eqn{rep-comod}{  \Hom_{\CV}(B^{\tens n},V)\to \Nat(F^n,F^n\tens V)}
This is the {\em representability assumption for comodules} and is always
satisfied
if $\CV$ is cocomplete and if the image of $F$ is rigid. In this case one can
write $B$ as a coend $B=\int^X F(X)^*\tens F(X)$.

\begin{theorem}\cite{Ma:bg} Let $F:\CC\to \CV$ obey the representability
assumption for comodules. Then $B$ is a bialgebra in $\CV$, uniquely
determined
as universal with the property that $F$ factors through $\CV^B$. If $\CC$ is
braided then $B$ is dual quasitriangular in the braided category with $\CR$
given by the ratio of the braidings in $\CC$ and $\CV$. If $\CC$ is rigid then
$B$ has a braided-antipode.
\end{theorem}
\proof Literally turn the above proofs up-side-down. \cite{Ma:bg} has more
traditional proofs. Note also the slightly different conventions there which
are chosen so as to ensure that the dual quasitriangular structure reduces for
$\CV={\rm Vec}$ to the usual notion as in Section~1.3 rather than its
convolution-inverse.
\endproof

By taking the identity functor to a cocompletion we obtain a canonical
braided-Hopf algebra
$B=\aut(\CC)$ associated to a rigid braided monoidal category
$\CC$\cite{Ma:bg}. By turning Corollary~3.10
up-side-down the braided-Hopf algebra this time is braided-commutative with
respect to a class of comodules. In this sense $\aut(\CC)$ is a braided group
of function algebra type.

\section{Applications to Ordinary Hopf Algebras}

In this section we give some applications of the above braided theory to
ordinary Hopf algebras. In this case there is either a background
quasitriangular bialgebra or Hopf algebra $H$ and we work in the braided
category ${}_H\CM$ or a background dual quasitriangular Hopf algebra $A$ and
we
work in the braided category $\CM^A$. See Sections~1.2 and~1.3 respectively.
The latter theory is the dual theory to the former and given by turning the
diagram-proofs in this book upside down. There need be no relationship between
$H$ or $A$ since we are dualizing the theory and not any specific Hopf
algebra.

In this context the content of Lemma~2.1 is immediate: it says that two
$H$-module algebras or $A$-comodule algebras have a braided tensor product
which is also an $H$-module algebra or $A$-comodule algebra respectively. From
(\ref{btens-stat}) and Theorem~1.10 or Theorem~1.16 we have obviously
\eqn{bratensH}{B,C,B\und\tens C\in {}_H\CM:\qquad (a\tens c)(b\tens d)=\sum
a(\CR\ut\la b)\tens (\CR\uo\la c)d}
\eqn{bratensA}{B,C,B\und\tens C\in \CM^A:\qquad(a\tens c)(b\tens d)=\sum a
b\bo\tens c\bo d\CR(c\bt\tens b\bt)}
for all $a,b\in B,\ c,d\in C$. We have introduced this construction in
\cite{Ma:bra}--\cite{Ma:exa} and explained there that it is exactly
a generalization of the $\Z_2$-graded or super-tensor product of
superalgebras.
This in turn has specific applications and spin-offs. An amusing one is

\begin{propos}\cite{BasMa:unb} Let $H$ be a Hopf algebra. Then the usual
$n$-fold tensor product $H^{\tens^n}$ is an $H$-module algebra under the
action
\[ h\la(b_1\tens \cdots\tens b_n)=\sum h_{(1)} b_1 Sh_{(2n)}\tens h_{(2)}b_2
Sh_{(2n-1)}\tens\cdots\tens h_{(n-1)}b_{n-1} Sh_{(n+2)}\tens h_{(n)}b_n
Sh_{(n+1)}.\]
\end{propos}
\proof This arises from the general theorems below as the $n$-fold braided
tensor product of $H$ as an $H$-module algebra and in the case where $H$ is
quasitriangular. This $n$-fold braided tensor product turns out to be
isomorphic to the usual tensor product in a non-trivial way using the
quasitriangular structure. Computing
the resulting action through this isomorphism and using the axioms
(\ref{qua2})
etc one finds the stated action. On the other hand the resulting formula does
not require any quasitriangular structure at all and can then be checked
directly to work for any Hopf algebra as stated.\endproof

This is a typical example in that the braided theory leads one to unexpected
formulas (as far as I know the last proposition is unexpected) which can then
be verified directly.  We will come to another such spin-off in Section~4.3.

\subsection{Transmutation}

In this section we shall show how to obtain non-trivial examples of bialgebras
and Hopf algebras in braided
categories. The key construction is one that we have called {\em
transmutation}
because it asserts that an
ordinary Hopf algebra can be turned by this process into a braided one. This
is achieved as an application of the  generalised reconstruction theorem in
Section~3 and is actually part of a rather general principle: by viewing an
algebraic structure in terms of its representations, and targeting these by
means of a functor into some new category, we can reconstruct our algebraic
object in this new category. In this way the category in which an algebraic
structure lives can be changed. Moreover, this `mathematical alchemy' can be
useful in that the structure may look more natural and have better properties
after transmutation to the new category.

In the present setting the data for our transmutation is a pair $H_1{\buildrel
f\over\to} H$ where $H_1$ is a quasitriangular Hopf algebra, $H$ is at least a
bialgebra, and $f$ a bialgebra map.

\begin{theorem}\cite{Ma:bra}\cite{Ma:tra} $H$ can be viewed equivalently as a
 bialgebra $B(H_1,H)$ living in the braided category ${}_{H_1}\CM$ by the
adjoint action induced by $f$. Here
\[ B(H_1,H)=\cases{H&{\rm as\ an\ algebra}\cr
 \und\Delta,\ \und S,\ \und\CR &{\rm modified\ coproduct,\ antipode,\
quasitriangular\ structure}}\]
where $B$ has a braided-antipode if $H$ has an antipode and a
braided-quasitriangular structure if $H$ is quasitriangular.
\end{theorem}
\proof We let $\CC={}_H\CM$ and $\CV={}_{H_1}\CM$ and $F$ the functor by
pull-back along $f$. Then Theorem~3.2 tells us there is a braided-Hopf algebra
$B$. \endproof

Explcit formulae for the transmuted structure are
\eqn{B(H)-hopf}{ \und\Delta b=\sum b\o f(S\CR_1\ut)\tens \CR_1\uo\la b\t,\quad
\und S b=\sum f(\CR_1\ut)S( \CR_1\uo\la b)}
\eqn{B(H)-univR}{ \und\CR=\sum \rho\uo f(S\CR_1\ut)\tens
\CR_1\uo\la\rho\ut,\qquad \rho= f(\CR_1^{-1})\CR.}
There is also an opposite coproduct characterised by
\eqn{B(H)-op}{\sum \Psi(b_{\und {(1)_{\rm op}}}\tens Q_1\uo\la
b_{\und {(2)_{\rm op}}})f(Q_1\ut)=\sum b_{\und {(1)}}\tens b_{\und {(2)}}}
where $Q_1=(\CR_1)_{21}(\CR_1)_{12}$ and $f(Q_1\ut)$ right-multiplies the
second tensor
factor of the output of $\Psi$. The underlines in the
superscripts
are to remind us that we intend here the braided-coproducts $\und\Delta$ and
$\und\Deltaop$. The
equation  can also be inverted to give an explicit formula for $\und\Deltaop$.
That these formulae obey the axioms (\ref{hopf-ax}) and (\ref{B-univR}) of
Section~2 is verified explicitly in \cite{Ma:tra}.

\begin{corol}\cite{Ma:any} Let $H$ be a Hopf algebra containing a group-like
element $g$ of order $n$. Then $H$ has a
corresponding anyonic version $B$. It has the same algebra and
\[\und\Delta b=\sum b\o g^{-\vert b\t\vert}\tens b\t,\quad
\und\eps b=\eps
b,\quad \und Sb=g^{\vert b\vert}Sb\]
\[\und\Delta^{\rm op} b=\sum b\t g^{-2\vert b\o\vert}\tens
g^{-\vert
b\t\vert}b\o,\quad \und\CR=\CR_{\Z_n'}^{-1}\sum\Ro g^{-\vert \Rt\vert}\tens
\Rt\]
\end{corol}
\proof We apply the transmutation theorem, Theorem~4.2 and compute the form of
$B=B(\Z_n',H)$. Here $\Z_n'$ is the non-standard quasitriangular Hopf algebra
in Example~1.7 with quasitriangular structure $\CR_{\Z_n'}$. The action of $g$
on $H$ is in the
adjoint representation $g\la b=gbg^{-1}$ for $b\in B$ and
defines the degree of homogeneous elements by $g\la b=q^{\vert
b\vert}b$.
\endproof

\begin{corol} Let $H$ be a quasitriangular Hopf algebra containing a
group-like
element $g$
of
order $2$. Then $H$ has a corresponding
super-version $B$.
\end{corol}
\proof The formulae are as in Corollary~4.3 with $n=2$. The commutation
relations with $g$ define the grading of $B$. \endproof

The first corollary was applied, for example to $H=u_q(g)$ at a root of unity
to simplify its structure. It led to a new simpler form for its
quasitriangular structure
 by finding its anyonic quasitriangular structure and working
back\cite{Ma:any}. The second corollary was usefully applied in
\cite{MaPla:uni} to superise the non-standard quantum group associated to the
Alexander-Conway polynomial. In these examples, a sub-quasitriangular Hopf
algebra is used to
generate the braided category in which the entire quasitriangular Hopf algebra
is then viewed by transmutation. In the process its quasitriangular structure
becomes reduced because the part from the
sub-Hopf algebra is divided out. This means that the part
corresponding to the sub-Hopf algebra is made in some sense cocommutative.

\begin{corol}\cite{Ma:bra} Every quasitriangular Hopf algebra $H$
has a braided-group analogue $B(H,H)$ which is braided-cocommutative in the
sense that $\und\CR=1\tens 1$ and $\und\Deltaop=\und\Delta$. The latter means
\[ \sum \Psi(b_{\und {(1)}}\tens Q\uo\la b_{\und
{(2)}})Q\ut=\sum b_{\und {(1)}}\tens b_{\und {(2)}}.\]
We call $B(H,H)$ the {\em braided group of enveloping algebra type} associated
to $H$. It is also denoted by $\und H$.
\end{corol}
\proof  Here we take the transmutation principle to its logical extreme and
view any quasitriangular Hopf algebra $H$ in its own braided category
${}_H\CM$, by
$H\subseteq H$. This is a bit like using a metric to determine geodesic
co-ordinates. In that co-ordinate system the metric looks locally linear.
Likewise, in its own category (as a braided group) our original
quasitriangular
Hopf algebra looks braided-cocommutative. From Corollary~3.10 we know that
$\Deltaop=\Delta$
and this gives the formula stated.
\endproof

This completely shifts then from one point of view (quantum=non-cocommutative
object in the usual category of vector spaces) to another
(classical=`cocommutative' but braided
object), and means that the theory of ordinary quasitriangular Hopf algebras
is
contained in the
theory of braided-groups.

The braided-Hopf algebra $B$ in Theorem~4.2 is equivalent to the original one
in that spaces
and algebras etc on which $H$ act also become transmuted to corresponding ones
for $B$. Partly, this is obvious since $B=H$ as an algebra, so any $H$-module
$V$ of $H$ is also a braided $B$-module. The key point is
that $V$ is also acted upon by $H_1$ through the mapping $H_1\to H$. So the
action of $H$ is used in two ways, both to define the corresponding action of
$B$ and to define the `grading' of $V$ as an object in a braided category
${}_{H_1}\CM$. This extends the process of
transmutation to modules.

\begin{propos}\cite[Prop~3.2]{Ma:bos} If $C$ is an $H$-module algebra then its
transmutation is
a braided $B(H,H)$-module algebra in the sense of (\ref{C-modalg}). Similarly
an $H$-module coalgebra becomes a braided $B(H,H)$-module-coalgebra. Here the
transmutation
does not change the action, but simply views it in the braided category.
\end{propos}
\proof An elementary computation from the form of (\ref{B(H)-hopf}) and the
braiding in Theorem~1.10. \endproof

For example, the adjoint action of $H$ on itself transmutes
to the braided-adjoint action of $B=B(H,H)$ on itself in Example~2.7.
Moreover, it means that $B(H,H)$ is braided-cocommutative with respect to
${\rm
Ad}$. Indeed

\begin{propos}\cite{Ma:bra}\cite{Ma:tra} For $B(H_1,H)$ the $\Deltaop$ behaves
like an opposite coproduct on all $B(H_1,H)$-modules that arise from
transmutation. In particular, $B(H,H)$ is cocommutative in the sense of
Definition~2.9 with respect to all braided-modules that arise
from transmutation.
\end{propos}
\proof Writing the braids in Definition~2.9 in terms of the quasitriangular
structure
 as explained in Theorem~1.10, we see that the condition for all $V$ is
implied by (and essentially equivalent to) the intrinsic
braided-cocommutativity formula in Corollary~4.5. \endproof

Finally, we mention a different aspect of this transmutation theory, namely a
result underlying the direct proof
that $B(H,H)$ is a braided-Hopf algebra (if one does not like the categorical
one).
\begin{lemma} Let $H$ be quasitriangular and $C$ be an algebra in ${}_H\CM$.
Then there is an algebra isomorphism
\[ \theta_{H,C}:H\tens C\to B(H,H)\und\tens C\]
where $\und\tens$ is the braided-tensor-product in Lemma~2.1.
\end{lemma}
\proof This is provided by $\theta_{H,C}(h\tens c)=\sum h S\CR\ut\tens
\CR\uo\la c$ and shows at once that $\und\Delta=\theta_{H,H}\circ\Delta$ in
(\ref{B(H)-hopf}) is an algebra homomorphism. The proof that $\theta$ is an
algebra
homomorphism is
\align{\theta_{H,C}(h\tens c)\theta_{H,C}(g\tens d)&=& \sum
h(S\CR\ut)(\CR''\ut\la(gS\CR'\ut))\tens  (\CR''\uo\CR\uo\la c)(\CR'\uo\la d)\\
&=&\sum h(S\CR\ut)\CR'''\ut g (S\CR'\ut)S\CR''\ut\tens
(\CR''\uo\CR'''\uo\CR\uo\la c)(\CR'\uo\la d)\\
&=&\sum hg S(\CR''\ut\CR'\ut)\tens(\CR''\uo\la c)(\CR'\uo\la
d)=\theta_{H,C}(hg\tens cd)}
using the definition (\ref{bratensH}) and the axioms (\ref{qua1}). \endproof

As an immediate example one has that $H^{\tens n}\isom B(H,H)^{\und\tens^n}$
for n-fold tensor products (iterate the lemma). Since $B(H,H)^{\und\tens^n}$
lives in ${}_H\CM$ as an algebra (an $H$-module algebra) via the adjoint
action, it follows that $H^{\tens^n}$ does also. Computing this gives
Proposition~4.1 as an amusing spin-off.

The transmutation theory obviously has a dual version for $A{\buildrel
f\over\to }A_1$ a bialgebra map
where $A_1$ is a dual quasitriangular Hopf algebra and $A$ is at least a
bialgebra.

\begin{theorem}\cite{Ma:eul}\cite{Ma:bg} $A$ can be viewed equivalently as a
 bialgebra $B(A,A_1)$ living in the braided category $\CM^{A_1}$ by the
right adjoint coaction induced by $f$. Here
\[ B(A,A_1)=\cases{A&{\rm as\ a\ coalgebra}\cr
 \und\cdot,\ \und S,\ \und\CR &{\rm modified\ product,\ antipode,\
dual quasitriangular\ structure}}\]
where $B$ has a braided-antipode if $A$ has an antipode and a braided-dual
quasitriangular structure if $A$ is dual quasitriangular.
\end{theorem}
\proof We let $\CC=\CM^A$ and $\CV=\CM^{A_1}$ and $F$ the functor by push-out
along $f$. Then Theorem~3.11 tells us there is a braided-Hopf algebra $B$. The
exact conventions for $\CR$ most useful here are in \cite{Ma:bg}.\endproof

Explicit formulae for the transmuted structure are\cite{Ma:bg}\cite{Ma:eul}
\eqn{B(A)-hopf}{a\und\cdot b=\sum a\t b\t \CR((S  a\o)a\th\tens S b\o),\quad
\und{S }a=\sum S  a\t\CR((S ^2a\th)S  a\o\tens a_{(4)})}
where for simplicity we concentrate on the case where $f$ is the identity (the
formulae in the general case are similar). The analogue of Corollary~4.5 is
that $B(A,A)$
is braided-commutative in the sense of Definition~2.9 turned up-side-down, for
all comodules $V$ that come from transmutation of comodules of $A$. This
reduces to an intrinsic form of commutativity dual to Corollary~4.5. This
comes
out explicitly as
\eqn{B(A)-com}{b\und \cdot a=\sum a\th\und\cdot b\th \CR(Sa\t\tens
b\o)\CR(a_{(4)}\tens b\t)\CR(b_{(5)}\tens Sa\o)\CR(b_{(4)}\tens a_{(5)}).}
We call $\und A=B(A,A)$ the {\em braided group of function algebra type
associated to $A$}. A direct proof that these formulae define a braided-Hopf
algebra as in (\ref{hopf-ax}) appears in \cite[Appendix]{Ma:bg}.

Finally, we suppose that $A$ is actually dual to $H$ in a suitable sense (for
example, one can
suppose they are finite dimensional). Until now we have not assumed anything
like this. Then from the two theorems above we have two Hopf algebras in
braided categories, and moreover the two categories can be identified in the
usual way.
Thus a right $A$-comodule defines a left $H$-comodule and viewing everything
in
this way in ${}_H\CM$ we have two Hopf
algebras $B(H,H)$ and $B(A,A)$ in the same braided category.

\begin{propos} If $A$ is dual to $H$ then the corresponding braided groups
$B(A,A)$ and $B(H,H)$ are dual in the braided category, $B(A,A)^\star=B(H,H)$.
\end{propos}
\proof  Explicitly, the duality is
given by $b\in B(H,H)$ mapping to a linear functional $<Sb,(\ )>$ on $B(A,A)$,
where $S$ is the usual antipode of $H$. See \cite{Ma:mec} for full details.
\endproof

Thus the usual duality if it exists becomes the categorical duality as in
Section~2.2. This is to be expected. More remarkable is the fact, also
verified
explicitly\cite{Ma:skl} that there is a canonical homomorphism of Hopf
algebras
in the braided category
\eqn{selfduality}{ Q:B(A,A)\to B(H,H),\qquad Q(a)=(a\tens
\id)(\CR_{21}\CR_{12})}
given by evaluation against $\CR_{21}\CR_{12}$. In the standard examples
$H=U_q(g)$ one has a formal expansion  $\CR_{21}\CR_{12}=1+2\hbar
K^{-1}+O(\hbar^2)$ where $K^{-1}$ is the inverse Killing form $g^*\to g$. So
(\ref{selfduality}) is a version for braided groups of this linear map
provided
by $Q$. This point of view has already been developed at the level of linear
maps $A\to H$ in \cite{ResSem:mat} where the Hopf algebra is called {\em
factorizable} if the map $Q$ is a linear isomorphism. What we have in
(\ref{selfduality}) is a much stronger statement: in the factorizable case
$B(A,A)\isom B(H,H)$ as Hopf algebras in a braided category. Since the first
of
these is the braided version of the quantum function algebra $\CO_q(G)$ and
the
second of the enveloping algebra $U_q(g)$, the isomorphism of their braided
versions is remarkable.

\subsection{Bosonization}

In this section we prove results going the other way, turning any Hopf algebra
in the braided category ${}_H\CM$ or $\CM^A$ into an ordinary Hopf algebra.
This process has been introduced by the author under the heading {\em
bosonization}. The origin of this term is from physics where $\Z_2$-graded
algebras etc are called `fermionic' while ordinary ungraded ones are called
`bosonic'. Not all
braided groups are of the type coming from the transmutation in the last
section, so bosonization is not simply transmutation in reverse.

\begin{theorem}\cite[Thm~4.1]{Ma:bos} Suppose that $B$ is a Hopf
algebra living in a braided category of the form ${}_H\CM$.
Then there is an ordinary Hopf algebra ${\rm bos}(B)=B\cocross H$.
\end{theorem}
\proof The abstract way that the result arose in \cite{Ma:bos} is as follows
(but once the
result is known a direct proof is also easy). Since $B$ lives in ${}_H\CM$ it
is in particular an $H$-module algebra. From Proposition~4.6 we see that the
same linear map makes $B$ a braided $B(H,H)$-module algebra. Likewise it is a
braided module
coalgebra and from
Proposition~4.7 this braided-module structure is braided-cocommutative. Hence
from Theorem~2.11 we have a semidirect product Hopf algebra $B\cocross
B(H,H)$.
This contains $B(H,H)$ and it is easy to see that this is indeed the
transmutation of $H\to H_2$ where $H_2$ is some ordinary Hopf algebra. It
computes explicitly as follows. As an algebra it is the semidirect or smash
product by the action of $H$ on $B$. So
$(1\tens h)(b\tens 1)=\sum h\o\la b\tens h\t$ where $\la$ is the action of
$H$.
As a coalgebra it is\cite{Ma:bos}
\eqn{boscoprod}{\Delta (b\tens h)=\sum b\Bo\tens \CR\ut h\o\tens \CR\uo\la
b\Bt\tens h\t.}
\endproof

Once the result and formula (\ref{boscoprod}) is known it is not hard to
verify
it directly. The key lemma for this direct verification is

\begin{lemma}\cite{Ma:dou} Let $H$ be a quasitriangular bialgebra or Hopf
algebra and $B$ a
left $H$-module with action $\la$. Then
\[ \beta(b)=\sum\CR\ut\tens \CR\uo\la b\]
makes $B$ into a left $H$-comodule. Moreover, it is compatible with $\la$ in
the sense of Example~1.3 and invertible so $B\in {}_H^H\CM$.
\end{lemma}
\proof Using the axioms of a quasitriangular bialgebra one sees at once that
this defines a coaction and this is compatible in the sense of Example~1.3. In
the case when $H$ is only
a bialgebra we have to check invertibility in the sense of (\ref{inv-comod}).
The required inverse is provided by ${\CR}^{-1}$ in place of $\CR$ in the
definition of $\beta$.
\endproof

As we explained in \cite{Ma:dou}, this defines a functor ${}_H\CM\to
{}_H^H\CM=\CZ({}_H\CM)$. Since the functor takes morphisms to morphisms (or by
direct computation) it is easy to see that if $B$ is an $H$-module (co)algebra
then it becomes in this way an $H$-comodule (co)algebra. It is completely
clear
then that ${\rm bos}(B)=B\cocross H$ has the structure of a semidirect
(co)product both as an algebra by $\la$, and as a coalgebra by the coaction
$\beta$ from  Lemma~4.12. This is the direct interpretation of
(\ref{boscoprod}). Simultaneous semidirect products and coproducts have been
studied in \cite{Rad:str} but the present construction of examples of them is
of course new and due to the author.

One thing that we learn from the categorical point of view is that this
ordinary Hopf algebra ${\rm bos}(B)$ is equivalent to the original $B$ in the
sense that its ordinary representations correspond to the
braided-representations of $B$\cite{Ma:bos}. This applies as much to
super-Hopf
algebras as to Hopf algebras on other categories, so we recover a result known
to experts working with super-Lie algebras and super-Hopf algebras that they
can be reduced to ordinary ones. See \cite{Fis:sch} for an example of this
strategy.

\begin{corol}\cite[Cor.~4.3]{Ma:bos} Any super-quasitriangular super-Hopf
algebra can be
bosonised to an
equivalent ordinary quasitriangular Hopf algebra. It consists of adjoining an
element $g$ with
relations $g^2=1$, $gb=(-1)^{|b|}bg$ and
\[\Delta g=g\tens g,\ \Delta b=\sum b\Bo g^{|b\Bt|}\tens b\Bt,\quad
S b=g^{-|b|}\und S b,\quad \CR=\CR_{\Z_2'}\sum \und\CR\uo
g^{|\und\CR\ut|}\tens
\und\CR\ut.\]
\end{corol}
\proof We have seen in Proposition~1.11 that the category of super-vector
spaces
is of the required form, with $H=\Z_2'$.
Here the $\Z_2$-graded modules of the original super-Hopf algebra are in
one-to-one correspondence with the usual representations of the bosonised
ordinary Hopf algebra. We suppose that we work over characteristic not $2$.
Also, we have written the formulae in a way that works in the anyonic case
with
2 replaced by $n$ and $(-1)$ by a primitive $n$-th root of unity for the Hopf
algebra structure.
\endproof

This means that the theory of super-Lie algebras (and likewise for colour-Lie
algebras\cite{Sch:gen}, anyonic quantum groups etc) is in a certain sense
redundant -- we could have worked with their bosonized ordinary Hopf algebras.
This is especially true in the super or colour case where there is no braiding
to complicate the picture.

As usual the above theory has exactly a version for a Hopf algebra living in
the braided category $\CM^A$ where $A$ is dual quasitriangular.
\begin{theorem} Suppose that $B$ is a Hopf
algebra living in a braided category of the form $\CM^A$.
Then there is an ordinary Hopf algebra ${\rm cobos}(B)=A\cross B$.
\end{theorem}
\proof We turn our diagram-proofs in the theory leading to Theorem~4.11
up-side-down. This time the coproduct is the semidirect one by the coaction
whereby $B$ is an object in $\CM^A$, and this also defines a right action
$\la$
(the dual version of Lemma~4.12) with respect to which we have a semidirect
product algebra on $A\tens B$,
\eqn{cobos}{b\ra a=\sum b\bo\, \CR(b\bt\tens a),\quad (1\tens b)(a\tens
1)=\sum
a\o\tens b\ra a\t.}
Just as one has a direct proof of Theorem~4.11, one can verify
directly that this ${\rm cobos}(B)$ is a Hopf algebra. \endproof

It is rather hard to consider this result and its attendant lemma as new
results since they are exactly the dual construction (by turning proofs as
diagrams up-side-down) of our bosonization Theorem~4.11. Equally well one
could
reflect in a mirror about a horizontal axis (or simply reverse the arrows in
the more conventional commutative diagrams). In this case left
modules/comodules become left comodules/modules and the Hopf algebra is
$\bar{{\rm cobos}}(B)=B\cocross A$ by a left handed semidirect product and
coproduct. Equally well
we could reflect in a vertical axis turning left modules to right modules etc
in Theorem~4.11 and giving a right-handed version $\bar{{\rm bos}}(B)=H\cross
B$.

As before there is no suggestion here that $A$ is dual to $H$ since it is the
construction that is being reversed and not any specific Hopf algebra. But if
$A$ is dual to $H$ (say finite-dimensional) then we can make both
constructions. Indeed, if $B^\star$ is the categorical dual as in Section~2.2
then
\eqn{cobosdual}{{\rm cobos}(B^\star)=A\cross B^\star\isom (B\cocross H)^*={\rm
bos}(B)^*.}
Details and (more importantly) an application may be found in \cite{Ma:mec}.
Another application of the bosonization theorem in its original form and in
the
dual form
can be found in \cite{CFW:sch} and \cite{FisMon:sch}, to prove nice
double-centraliser theorems for various kinds of Lie algebras in symmetric
categories.

\subsection{Radford's Theorem}

In \cite[Appendix]{Ma:skl} we have explained in detail how the above ideas
provide a new braided interpretation of Radford's theorem\cite{Rad:str} about
Hopf algebras with projections. This theorem asserts that if $H,H_1$ are
ordinary Hopf algebras and if
\eqn{proj}{H_1{
p\atop{{\longrightarrow\atop \hookleftarrow}\atop i}}H}
are bialgebra maps with $p\circ i=\id$ (a Hopf algebra projection), then there
is an algebra and
coalgebra $B$ such that $H_1\isom B\cocross H$ as a simultaneous semidirect
product and semidirect coproduct. Radford called such simultaneous semidirect
(co)products where the result is a Hopf algebra `biproducts' and showed that
they correspond to projections. We have already introduced some examples in
the
last section (arising from the bosonization process) but now we consider the
general situation.

At the time of \cite{Rad:str} the notion of braided categories
was yet to be invented. Because of this the algebra and coalgebra $B$ in
Radford's theorem was simply some exotic object where the algebra and
coalgebra
did not form an ordinary Hopf algebra. We have pointed out in \cite{Ma:skl}
and
cf.\cite{Ma:dou}
that $B$ is in fact nothing other than a Hopf algebra in the braided category
${}_H^H\CM={}_{D(H)}\CM$ of Example~1.3 (we stressed the latter in
\cite{Ma:skl} for pedagogical reasons but explained that the former was
more useful in the infinite-dimensional case). Thus we arrived at the
following interpretation of Radford's theorem.

\begin{propos} \cite[Prop. A.2]{Ma:skl} Let $H_1{
p\atop{{\longrightarrow\atop \hookleftarrow}\atop i}}H$ be a Hopf algebra
projection and let $H$ have invertible antipode. Then there is a Hopf algebra
$B$ living in the braided category
${}_H^H\CM$ such that $B\cocross H\isom H_1$.
\end{propos}
\proof Explicitly, $B$ is a subalgebra of $H_1$ and in ${}_H^H\CM$ by action
$\la$ and coaction $\beta$,
\eqn{RadBob}{ B=\{b\in H_1\ |\ \sum b\o\tens p(b\t)=b\tens 1\},\quad h\la
b=\sum i(h\o)b\ant \circ i(h\t),\quad \beta(b)=p(b\o)\tens b\t}
where $h\in H$. The braided-coproduct, braided-antipode and braiding of $B$
are
\eqn{RadB-hopf}{\und\Delta b=\sum b\o\ant \circ i\circ p(b\t)\tens b\th,\quad
\und\ant b=\sum i\circ p(b\o)\ant b\t,\quad\Psi_{B,B}(b\tens c)=\sum p(b\o)\la
c\tens b\t.}
The isomorphism $\theta:B\cocross H\to H_1$ is $\theta(b\tens h)=b i(h)$, with
inverse $\theta^{-1}(a)=\sum a\o\ant\circ i\circ p(a\t)\tens p(a\th)$ for
$a\in
H_1$. The only new part beyond \cite{Rad:str} is the identification of the
`twisted Hopf algebra' $B$ now as a Hopf algebra living in a braided category,
and some slightly more explicit formulae for its structure. The set $B$
coincides with the image of the projection $\Pi:H_1\to H_1$ defined by
$\Pi(a)=\sum a\o \ant \circ i\circ p(a\t)$ in \cite{Rad:str}, while the
pushed-out left adjoint coaction of $H$ on $B$ then reduces to the left
coaction as stated. The braiding is from Example~1.3. The axioms of a Hopf
algebra in a braided category require that $\und \Delta:B\to B\tens B$ is an
algebra homomorphism with respect to the braided tensor product algebra
structure on $B\tens B$. Writing $\und\Delta b=\sum b_{\und{(1)}}\tens
b_{\und{(2)}}$, this reads
\eqn{rad-hopf}{ \und\Delta (bc)=\sum b_{\und{(1)}}\Psi(b_{\und{(2)}}\tens
c_{\und{(1)}})c_{\und{(2)}}=\sum b_{\und{(1)}}\, (b_{\und{(2)}}\bo\la
c_{\und{(1)}})\tens b_{\und{(2)}}\bt c_{\und{(2)}}}
which indeed derives the condition in \cite{Rad:str}. The structure of
$B\cocross H$ is the standard left-handed semidirect one by the action and
coaction stated. Applying $\theta$ to these structures and evaluating further
at once gives $\theta$ as a Hopf algebra isomorphism. Of course if $H_1$ is
only a bialgebra then $B$ is only a bialgebra in a braided category. In this
case one can
use the convolution inverse $i\circ S$ in the above. Also, the restriction to
invertible antipode on $H$ is needed only to ensure that $\Psi$ is invertible
as explained in Example~1.3. It is part of our interpretation of $B$ as a
braided-Hopf algebra rather than part of Radford's theorem itself. \endproof

Such Hopf algebra projections have a geometrical interpretation as examples of
trivial quantum principal bundles  \cite{BrzMa:gau} and at the same time as
quantum mechanics\cite{Ma:mec}. These papers also make some limited contact
with more established ideas of non-commututive geometry as in \cite{Con:alg}.
Also note that in the above the Hopf algebra $H$ need not be quasitriangular
or
dual quasitriangular. If it is then the above construction becomes related to
the bosonization theorems of the last subsection, as mentioned there.

\section{Braided Linear Algebra}

In this section we describe some general constructions for examples of
bialgebras and Hopf algebras in braided categories associated to
a general matrix solution of the QYBE as in Example~1.1. This includes some
interesting algebras for ring theorists. The first of these, $B(R)$, has a
matrix of generators
with matrix coproduct and includes a degenerate form of the Sklyanin algebra
for the usual $GL_q(2)$ $R$-matrix. The second has a vector or covector of
generators with linear coproduct, and includes the famous quantum plane
$yx=qxy$ for this $R$-matrix. Thus the quantum plane does have a linear
addition law provided we work in a braided category. Finally, we mention some
recent developments such as a notion of braided-Lie algebra.

\subsection{Braided Matrices}

Let $R$ be an invertible matrix solution of the QYBE and $A(R)$ the associated
dual quasitriangular bialgebra as in Section~1.3. In the nicest case we can
quotient $A(R)$ to obtain a dual quasitriangular Hopf algebra $A$. We say in
this case that $R$ is {\em regular}. This is true for the standard $R$
matrices
and one obtains $A=\CO_q(G)$ as shown in \cite{FRT:lie}. Another way to obtain
a dual quasitriangular Hopf algebra is if $R$ is dualizable. In this case we
have an
associated dual quasitriangular Hopf algebra $A=GL(R)$. Either way we have a
canonical bialgebra map $A(R)\to A$ and can apply the transmutation theorem,
Theorem~4.9
to obtain a bialgebra $B(R)=B(A(R),A)$ in the category $\CM^A$. This gives the
following construction, which we verify directly.

\begin{propos}\cite{Ma:eul}\cite{Ma:exa} Let $R$ be a bi-invertible solution
of
the QYBE with $v^i{}_j=\widetilde{R}{}^i{}_a{}^a{}_j$ also invertible. Then
there is a bialgebra $B(R)$ in the
braided category of $A$-comodules, with matrix generators $\vecu=\{u^i{}_j\}$
and relations, braiding and coalgebra
\[ R^k{}_a{}^i{}_b u^b{}_c R^c{}_j{}^a{}_d u^d{}_l=u^k{}_a
R^a{}_b{}^i{}_c u^c{}_d R^d{}_j{}^b{}_l,\quad {\rm i.e.}\quad
R_{21}\vecu_1R_{12}\vecu_2= \vecu_2 R_{21} \vecu_1 R_{12}.\]
\[\Psi(u^i{}_j\tens u^k{}_l)=u^p{}_q\tens u^m{}_n  R^{i}{}_a{}^d{}_{p}
R^{-1}{}^a{}_{m}{}^{q}{}_b
R^{n}{}_c{}^b{}_{l} {\tilde R}^c{}_{j}{}^{k}{}_d,\quad {\rm i.e.}\quad
\Psi(R^{-1}\vecu_1 \tens R\vecu_2)=\vecu_2 R^{-1}\tens \vecu_1 R\]
\[ \und\Delta u^i{}_j=u^i{}_a\tens u^a{}_j,\ \und\eps
u^i{}_j=\delta^i{}_j\quad{\rm i.e.}\quad \und\Delta\vecu=\vecu\tens\vecu,\quad
\und\eps\vecu=\id \]
\end{propos}
\proof The formulae are obtained from Theorem~4.9 and then verified directly.
Bi-invertible means $R^{-1}$ and the second-inverse $\tilde{R}$ exist and is
all that is
needed to verify the braided-bialgebra axioms. The existence of $v^{-1}$ is
needed in defining
$\Psi^{-1}$ and is equivalent to demanding that $R$ is dualizable.
The commutation relations come from (\ref{B(A)-com}) using Proposition~1.15 to
evaluate the formulae in matrix form. This gives a matrix equation of the
form $uu=uu R^{-1}RR\widetilde R$ with appropriate indices, see
\cite{Ma:eul}\cite{Ma:exa}.
Putting two of the
$R$'s to the left (or rearranging (\ref{B(A)-com})) gives the formula stated.
The same applies to the braiding which comes out as a product of 4
$R$-matrices
as shown but is also conveniently written in the compact form stated. Its
extension to products is by definition in such a way that the product is a
morphism in the category generated by this braiding $\Psi$ and one can
verify that this extension is consistent with the relations. To verify that
the
result is indeed a braided bialgebra an even more compact notation is
useful. Namely as in (\ref{btens-stat}) we can label the second copy of $B(R)$
in the braided tensor product by a prime, and then suppress $\Psi$. Then the
commutation relations in the braided tensor product $B(R)\und\tens B(R)$ are
$R^{-1}\vecu'_1 R\vecu_2=\vecu_2 R^{-1}\vecu'_1 R$ and we compute
\align{&&\nquad R_{21}\vecu_1\vecu'_1
R\vecu_2\vecu'_2=R_{21}\vecu_1 R(R^{-1}\vecu'_1
R\vecu_2)\vecu'_2=(R_{21}\vecu_1 R\vecu_2) R^{-1}R_{21}^{-1}(R_{21}\vecu'_1
R\vecu'_2)\\
&&= \vecu_2 R_{21}(\vecu_1 R_{21}^{-1}\vecu'_2 R_{21})\vecu'_1
R=\vecu_2 R_{21}R_{21}^{-1}\vecu'_2 R_{21}\vecu_1\vecu'_1 R
= \vecu_2\vecu'_2 R_{21}\vecu_1\vecu'_1 R}
as required for $\und\Delta$ to extend to $B(R)$ as a bialgebra in a braided
category. In each expression, the brackets indicate how to apply the relevant
relation to obtain the next expression. \endproof

Note that the braided category here is generated by the matrix in
$M_{n^2}\tens M_{n^2}$ corresponding to $\Psi$ as stated. On the other hand in
the regular or dualizable case we can identify this category as contained in
that of $A$-comodules, under the induced adjoint coaction in Theorem~4.9. In
the present setting this is
\eqn{Ad-coact-mat}{\beta(u^i{}_j)=u^m{}_n \tens (St^i{}_m)t^n{}_j,\qquad{\rm
i.e.}\quad \vecu\to\vect^{-1}\vecu\vect.}

Note also that if $R_{21}R_{12}=1$ (the triangular case) then the braiding and
the commutation relations co-incide so that $\cdot=\Psi^{-1}\circ\cdot$ which
is the naive notion of braided-commutativity. In general however this will not
do and instead the braided-commutativity relations  are different from the
braiding itself.

Finally, it is good to know that the  algebra $B(R)$ has at least one
canonical
representation. This is provided by
\eqn{B(R)-can}{\rho(u^i{}_j)^k{}_l=Q^i{}_j{}^k{}_l,\quad {\rm i.e.}\quad
\rho_2(\vecu_1)=Q_{12};\quad Q=R_{21}R_{12}}
In the compact notation the proof reads
\align{\rho_3(R_{21}\vecu_1R_{12}\vecu_2)&=&R_{21}
\rho_3(\vecu_1)R_{12}\rho_3(\vecu_2)\\
&=& R_{21} Q_{13} R_{12} Q_{23}= Q_{23} R_{21} Q_{13} R_{12}\\
&=& \rho_{3}(\vecu_2)R_{21} \rho_3(\vecu_1)R_{12}=\rho_3(\vecu_2 R_{21}\vecu_1
R_{12}).}
The middle equality follows from repeated use of the QYBE. Note that this
representation is trivial in the triangular case. Other useful representations
can be built from this.

{}From the theory above we obtain also that these braided matrices are related
to usual
quantum matrices by transmutation. The formula for the modified product comes
out from (\ref{B(A)-hopf}) as \cite{Ma:lin}
\eqn{u-t-trans}{\begin{array}{rcl}u^i{}_j&=&t^i{}_j\\  u^i{}_j
u^k{}_l&=&t^a{}_bt^d{}_lR^i{}_a{}^c{}_d \widetilde R^b{}_j{}^k{}_c \\ u^i{}_j
u^k{}_l u^m{}_n&=&  t^d{}_b t^s{}_u t^z{}_n R^i{}_a{}^p{}_q R^a{}_d{}^w{}_y
{\widetilde
R}^b{}_c{}^v{}_w {\widetilde R}^c{}_j{}^k{}_p R^q{}_s{}^y{}_z
{\widetilde R}^u{}_l{}^m{}_v \end{array}}
etc. Here the products on the left are in $B(R)$ and are related by
transmutation to the products on the right which are in $A(R)$. If we write
some or all of the $R$-matrices over to the left hand side we have equally
well
the compact matrix form \cite{Ma:lin},
\eqn{t-u-trans}{\begin{array}{rcl}\vecu &=&\vect\\ R^{-1}_{12}\vecu_1
R_{12}\vecu_2& = &\vect_1\vect_2 \\
R_{23}^{-1} R_{13}^{-1}R_{12}^{-1}\vecu_1 R_{12}\vecu_2
R_{13}R_{23}\vecu_3&=&
\vect_1\vect_2\vect_3 \end{array}}
etc. This is just a rearrangement of (\ref{u-t-trans}) or our universal
formula
(\ref{B(A)-hopf}). For the transmuted product of multiple strings the
universal
formula from (\ref{B(A)-hopf}) involves a kind of partition function made from
products of $R$ to transmute the bosonic $A(R)$ to the braided $B(R)$
\cite{Ma:lin}.

We believe these braided matrices deserve more study as certain well-behaved
quadratic algebras. For the standard $R$-matrices they have quotients giving
the braided versions $B_q(G)$ say of the quantum function algebras $\CO_q(G)$,
and are at the same time isomorphic via (\ref{selfduality}) and for generic
$q$
to the braided versions of the corresponding quantum enveloping algebras
$U_q(g)$. This is related to constructions in physics\cite{ResSem:cen}. On the
other hand these $B(R)$ are interesting even at the quadratic level. The case
of $BM_q(2)$ for the $GL_q(2)$ $R$-matrix was studied in \cite{Ma:exa} and
shown in \cite{Ma:skl} to be a degenerate form of the Sklyanin algebra. Some
remarkable homological properties of these braided matrices have recently been
obtained in \cite{LeB:hom}.

\subsection{Braided Planes}

If the algebras $B(R)$ are like `braided matrices' because they have a matrix
of generators with matrix coproduct, one can complete the picture with some
notion of `braided vectors' and `braided covectors'. The usual algebra
suggested here from physics is the Zamalodchikov or `exchange algebra'  of the
form $\vecx_1\vecx_2=\lambda\vecx_2\vecx_1R$ where $\vecx$ say is a row vector
of generators and $R$ is a matrix obeying the QYBE, and $\lambda$ is a
normalization constant. This is an interesting algebra in the Hecke
case\cite{Gur:alg} but does not seem so interesting for a general solution of
the QYBE. Coming out of the theory of Hopf algebras in braided categories it
would seem rather more natural to use $R$ for the braiding and some slightly
different matrix $R'$ say for the commutation relations of the algebra. This
strategy works and gives the following
construction.

\begin{propos}\cite{Ma:poi} Let $R\in M_n\tens M_n$ be an invertible solution
of the QYBE and $R'\in M_n\tens M_n$ an invertible  matrix obeying
\[ R_{12}R_{13}R'_{23}=R'_{23}R_{13}R_{12},\quad
R_{23}R_{13}R'_{12}=R'_{12}R_{13}R_{23},\quad (PR+1)(PR'-1)=0\]
 where $P$ is the permutation matrix.
Then there are braided-bialgebras $\Vhaj(R')$ and $V(R')$ with row and
column vectors of generators $\vecx=\{x_i\}$ and $\vecv=\{v^i\}$ respectively
and with relations and braiding
\[ x_ix_j=x_bx_a R'^a{}_i{}^b{}_j,\  \Psi(x_i\tens x_j)=x_b\tens x_a
R^a{}_i{}^b{}_j,\quad {\rm i.e.}\
 \vecx_1\vecx_2=\vecx_2\vecx_1R',\ \Psi(\vecx_1\tens\vecx_2)=\vecx_2\tens
\vecx_1R\]
\[ v^iv^j=R'^i{}_a{}^j{}_bv^bv^a,\  \Psi(v^i\tens v^j)=R^i{}_a{}^j{}_b
v^b\tens
v^a,\quad {\rm i.e.}\
 \vecv_1\vecv_2=R'\vecv_2\vecv_1,\ \Psi(\vecv_1\tens\vecv_2)=R\vecv_2\tens
\vecv_1\]
and linear coalgebra
\[\und\Delta \vecx=\vecx\tens 1+1\tens
\vecx,\quad\und\eps\vecx=0,\quad\und\Delta\vecv=\vecv\tens
1+1\tens\vecv,\quad\und\eps\vecv=0.\]
If $R_{21}R'_{12}=R'_{21}R_{12}$ then these are braided-Hopf algebras with
$\und S \vecx=-\vecx$ and $\und S\vecv=-\vecv$
\end{propos}
\proof One can see that $\Psi$ extends to products of generators in such a way
that the product is a morphism in the braided category generated by $R$. For
details see \cite{Ma:poi}. To see that the result forms a bialgebra we use the
compact notation as in the proof of Proposition~5.1,
$(\vecx_1+\vecx'_1)(\vecx_2+\vecx'_2)=\vecx_1\vecx_2+\vecx'_1\vecx_2+\vecx_1
\vecx'_2+\vecx'_1\vecx'_2=\vecx_2\vecx_1R'_{12}+\vecx_1\vecx'_2 (PR_{12}+1)
+\vecx'_2\vecx'_1R'_{12}$ while $(\vecx_2+\vecx'_2)(\vecx_1+\vecx'_1)R'_{12}$
has the same outer terms and the cross terms $\vecx'_2\vecx_1R'_{12}+\vecx_2
\vecx'_1 R'_{12}=\vecx_1\vecx'_2(R_{21}+P)R'$. These are equal since
$PR+1=PRPR'+PR'$ from the assumption on $R'$. Similarly for the other details.
\endproof

There are lots of ways to satisfy the auxiliary equation for the matrix $R'$.
If $R$ is a Hecke symmetry we can simply take $R'=\lambda R$ where $\lambda$
is
a suitable normalization. This is the familiar case. For example the standard
$GL_q(n)$ $R$-matrix gives the algebra
\eqn{qplane}{ x_ix_j=qx_jx_i \quad {\rm if}\quad i>j}
which is the n-dimensional quantum plane. Thus we see that it forms a Hopf
algebra in the braided category generated by $R$. There is another
normalization $\lambda$ giving another quantum plane algebra and again a
braided-Hopf algebra.

Another solution is $R'=P$ the permutation matrix. In this case the relations
of the algebra are free (no relations). The $R$-matrix still enters into the
braiding. This free-braided plane is therefore canonically associated to any
invertible matrix solution of the QYBE. The simplest member of this family is
the {\em braided line}. This is $B=k<x>$ (one generator) and
\eqn{braline}{\Psi(x\tens x)=qx\tens x,\quad \und\Delta x=x\tens 1+1\tens
x,\quad \und\eps x=0,\quad \und S x=-x.}

Finally, any $R$ will obey an equation of the form $\prod_i(PR-\lambda_i)=0$
and for each $\lambda_i$ we can rescale $R$ so that one of these factors
become
$PR+1$. Then $PR'-1$ defined as a multiple of the remaining factors gives us
an
$R'$ obeying the required equations. Thus there are canonical braided Hopf
algebras
$\Vhaj(R')$
and $V(R')$ for  each eigenvalue $\lambda_i$ in the decomposition. In the
Hecke
case there are by definition two such eigenvalues corresponding to the two
normalizations mentioned, but in general $R'$ will not be a multiple of $R$
and
need not obey the QYBE.

The algebra $V(R')$ can be called the {\em braided vector} algebra. Likewise
the algebra $\Vhaj(R')$ can be called the {\em braided  covector} algebra.
Another notation is $V^*(R')$ but we have avoided this here in order to not
suggest that it is the left dual of the braided vector algebra in the sense of
Section~2.2. In fact, the braiding is such that in the dualizable case the
vector space generating the vector algebra is the left dual of that generating
the covector algebra (compare with Section~1.1). Together with the braided
matrices $B(R)$, these algebras form a kind of braided linear algebra. We
refer
to \cite{Ma:lin} for more details.

\subsection{New Directions}

In the above we have described the basic theory of Hopf algebras in braided
categories and some canonical canonical constructions for them which have been
established for the most part in the period 1989 -- 1991. Some more recent
directions are as follows.

One direction is that of {\em braided differential calculus}. The idea is that
just as the differential structure on a group is obtained by making an
infinitesimal group translation, now we can use our braided-coproducts on the
braided-planes  and braided-matrices to obtain corresponding differential
operators.

\begin{propos}\cite{Ma:fre} The operators $\del^i:\Vhaj(R')\to \Vhaj(R')$
defined by
\[ \del^ix_{i_1}\cdots x_{i_m}=\delta^i{}_{j_1}x_{j_2}\cdots
x_{j_m}\left[m;R\right]^{j_1\cdots j_m}_{i_1\cdots i_m}\]
\[\left[m;R\right]=1+(PR)_{12}+(PR)_{12}(PR)_{23}
+\cdots+(PR)_{12}\cdots (PR)_{m-1,m}\]
obey the relations of $V(R')$ and the braided-Leibniz rule
\[ \del^i(ab)=(\del^i a)b+\cdot\Psi^{-1}(\del^i\tens a)b,\qquad \forall a,b\in
\Vhaj(R').\]
\end{propos}

The result says that the braided-covector algebra (which is like the algebra
of
co-ordinate functions on some kind of non-commutative algebraic variety) is a
left $V(R')$-module algebra in the category with reversed braiding.
For the braided-line we obtain the usual Jackson $q$-derivative $(\del
f)(x)={f(qx)-f(x)\over(q-1)x}$. For the famous quantum plane we recover its
well-known two dimensional differential calculus obtained usually by other
means. We have also introduced in this context a braided-binomial theorem for
the `counting' of braided partitions.  This is achieved by means of the {\em
braided integers} $\left[m,R\right]$ in the proposition. One can also define a
braided exponential map $\exp_R(\vecx|\vecv)$ and prove a braided Taylor's
theorem at least in the free case where $R'=P$. We refer to \cite{Ma:fre} for
details.

Likewise, one has some natural right-handed differential operators on the
braided
matrices obtained from Example~2.8 applied to $B(R)$ and its braided-Hopf
algebra quotients (the result lifts to the bialgebra setting). These are
computed
in detail in \cite{Ma:lie}.

Related to this, one can formalise at least one notion of {\em braided Lie
algebra}\cite{Ma:lie}. In a symmetric (not braided) category the notion of a
Lie
algebra and its enveloping algebra are just the obvious ones,
see\cite{Gur:yan}.
But as soon as one tries to make this work in the
braided setting there are problems with the naive approach. One solution is
based on the notion of braided-cocommutativity in Definition~2.9 with respect
to the braided-adjoint action. Based on its properties in this case one can
extract the axioms of a braided-Lie algebra as the following: a coalgebra
$(\CL,\Delta,\eps)$ in the braided category and a bracket $[\ ,\
]:\CL\tens\CL\to\CL$ obeying
\eqn{lie-ax}{\epsfbox{Lie.eps}}
We show that such a braided-Lie algebra has an enveloping bialgebra living in
the category.

The braided matrices $B(R)$ (e.g. the degenerate Sklyanin algebra) double up
in
this role as braided-enveloping bialgebras. Here $\CL=\span\{u^i{}_j\}$ with
the matrix coproduct as in Section~5.1. These generators are a mixture of
`group-like' and `primitive-like' generators. If one wants something more
classical one can work equally well with the
space $\CX=\span\{\chi^i{}_j\}$ where $\chi^i{}_j=u^i{}_j-\delta^i{}_j$. In
these terms the braiding has the same form as between the $\vecu$ generators,
while the relations and coproduct become
\eqn{chi-lie}{R_{21}\chi_1R_{12}\chi_2-\chi_2 R_{21}\chi_1
R_{12}=Q_{12}\chi_2-\chi_2Q_{12},\quad\und\Delta\chi=\chi\tens
1+1\tens\chi+\chi\tens\chi}
in the compact notation. This is more like a `Lie algebra' and for the
standard
$R$ matrices, as the deformation parameter $q\to 1$ the matrix
$Q=R_{21}R_{12}\to \id$ and the commutator on the right hand side vanishes.
This means that it is the
rescaled generators $\bar\chi=\hbar^{-1}\chi$ that tend to a usual Lie algebra
and
indeed the extra non-primitive term in the coproduct of the $\bar\chi$ now
tends to
zero. For details see \cite{Ma:lie}. We recall also that for the standard $R$
matrices the algebras $B(R)$ have a quotient which is $U_q(g)$ at least
formally, so these are understood as (the quotient of) the braided-enveloping
algebra of a braided-Lie algebra.

We have not tried here to survey a number of more specific applications of
this
braided work. These
include \cite{BrzMa:cla}\cite{GurMa:bra}\cite{LyuMa:bra}\cite{Ma:inf} as well
as
more physically-based applications. In addition among many relevant and
interesting
recent works by other mathematicians and that I have not had a chance to touch
upon, I would like
to at least mention
\cite{Fre:hig}\cite{KapVoe:bra}\cite{Mor:som}\cite{ShnSte:cob} in the
categorical
direction and \cite{CohWes:sup}\cite{Rad:sol} as well as
\cite{BMO:mul}\cite{CFW:sch}\cite{FisMon:sch}\cite{LeB:hom} already mentioned
for works in an algebraic direction.
What we have shown is that a number of
general mathematical constructions can be braided. There is clearly plenty of
scope for further work in this programme. Some potential areas are
applications to knot theory, a general braided-combinatorics (based on the
braided-integers above) and some kind of braided-analysis.


\end{document}